\documentclass[preprint]{elsarticle}
\usepackage{amssymb}
\usepackage{amsmath,amsthm}
\usepackage{amstext}
\usepackage{mathptmx}
\usepackage{bm}
\usepackage{booktabs}
\usepackage{multirow, mathtools}
\usepackage{lineno,hyperref}
\usepackage{color}
\usepackage{flafter}
\usepackage{tabularx}
\usepackage[hang]{subfigure}
\usepackage{algorithm}
\usepackage{algpseudocode}
\modulolinenumbers[5]

\newtheorem{theorem}{Theorem}

\usepackage{pgfplots}
\pgfplotsset{compat=newest}
\usetikzlibrary{plotmarks}
\usetikzlibrary{arrows.meta}
\usepgfplotslibrary{patchplots}
\usepackage{grffile}
\newlength\fwidth

\journal{a journal}

\graphicspath{{figures/}}

\begin{document}

\begin{frontmatter}
\title{Tensor train-Karhunen-Lo{\`e}ve expansion for continuous-indexed random fields using higher-order cumulant functions}


\author[mymainaddress]{Ling-Ze Bu\corref{mycorrespondingauthor}}
\cortext[mycorrespondingauthor]{Corresponding author}
\ead{17b933010@stu.hit.edu.cn}
\author[mymainaddress,mysecondaryaddress,mythirdaddress]{Wei Zhao}
\author[mymainaddress,mysecondaryaddress,mythirdaddress]{Wei Wang}
\ead{wwang@hit.edu.cn}
\address[mymainaddress]{School of Civil Engineering, Harbin Institute of Technology, Harbin 150090, China}
\address[mysecondaryaddress]{Key Lab of Structures Dynamic Behavior and Control of the Ministry of Education, Harbin Institute of Technology, Harbin 150090, China}
\address[mythirdaddress]{Key Lab of Smart Prevention and Mitigation of Civil Engineering Disasters of the Ministry of Industry and Information Technology, Harbin Institute of Technology, Harbin 150090, China}

\begin{abstract}
	
The goals of this work are two-fold: firstly, to propose a new theoretical framework for representing random fields on a large class of multidimensional geometrical domain in the tensor train format; secondly, to develop a new algorithm framework for accurately computing the modes and the second and third-order cumulant tensors within moderate time. The core of the new theoretical framework is the tensor train decomposition of cumulant functions. This decomposition is accurately computed with a novel rank-revealing algorithm. Compared with existing Galerkin-type and collocation-type methods, the proposed computational procedure totally removes the need of selecting the basis functions or collocation points and the quadrature points, which not only greatly enhances adaptivity, but also avoids solving large-scale eigenvalue problems. Moreover, by computing with third-order cumulant functions, the new theoretical and algorithm frameworks show great potential for representing general non-Gaussian non-homogeneous random fields. Three numerical examples, including a three-dimensional random field discretization problem, illustrate the efficiency and accuracy of the proposed algorithm framework.

\end{abstract}


\end{frontmatter}


\section{Introduction}\label{Intro}

Uncertainty quantification in engineering and applied sciences often requires random field descriptions of spatial variability of uncertain media and loads, etc. Essentially, random field representation is a problem of \emph{data compression and reconstruction}. Information of a random field must be compressed in a sufficiently small number of deterministic functions and random variables which can reconstruct the random field with acceptable accuracy. 

A frequently used approach is the Karhunen-Lo{\`e}ve (K-L) expansion \cite{Karhunen1947,Loeve1948} (also named as Principal Component Analysis (PCA) in the statistical community) which was originally proposed for representing stochastic processes. A continuous-time stochastic process is represented with countable modes and latent factors. To model spatial variability of a physical quantity in a multidimensional space, the most common approach is directly extending K-L expansion by treating the spatial coordinate as a macro parameter (named as trivial PCA throughout this paper). The modes are eigenfunctions of an integral operator induced by the covariance function and the latent factors are uncorrelated in the sense of second-order moments. A comprehensive review of the methods (including Galerkin-type and collocation-type ones) for computing the modes was made in \cite{Betz2014}. Other procedures such as spectral element projection \cite{Oliveira2014}, multilevel finite element method \cite{Xie2015}, two-dimensional Haar wavelet Galerkin method \cite{Azevedo2018} and high-order polynomial-based Ritz-Galerkin method \cite{Zhang2019} were also proposed in recent years. To tackle the difficulties of the previous methods in treating complex geometries, modes were projected onto a subspace spanned by isogeometric basis functions in recent papers \cite{Rahman2018,Li2018,Jahanbin2019}. Meanwhile, several recent works tried to avoid the computationally expensive eigenvalue problem in trivial PCA based on the idea of variable separation. Ghosh et al. \cite{Ghosh2014} redefined the covariance function of a spatial-temporal random field, and computed both the temporal and spatial modes with CP and Tucker decomposition, respectively. Zentner et al. \cite{Zentner2016} and Guo et al. \cite{Guo2016} independently proposed a hierarchical orthogonal decomposition method to split the temporal and spatial modes. Similar idea was also shown in \cite{Zheng2017}, however, this work did not really achieve variable separation since the modes still contain multiple coordinates. The stepwise covariance matrix decomposition method \cite{Li2019} derives a Tucker-type representation and is only suitable for rank-1 covariance functions. 

Latent factors also play an important role for describing the probabilistic structure of the random field. For non-Gaussian random fields which widely exist in practice, the latent factors are generally non-Gaussian and have higher-order dependencies, making the representation of their probabilistic structure a nontrivial task. To achieve this task, several methods are dedicated to match the prescribed single-point marginal distribution function and covariance function. Phoon et al. \cite{Phoon2005} and Dai et al. \cite{Dai2019} proposed similar iterative procedures to update the marginal distribution of each latent factor. Kim et al. \cite{Kim2015a} iteratively updated the covariance function of the underlying Gaussian process based on the translation process theory. Other works tried to estimate the probabilistic structure of the latent factors with data-driven procedures \cite{Ganapathysubramanian2008,Ma2011,Poirion2014,Olivier2016}. Since this paper is focused on representing a random field with prescribed statistics, we will not give further comments on these methods. 

For the computation of modes, although the complex physical domains can be handled by the aforementioned isogeometric analysis (IGA)-based methods, however, the underlying architecture is still the trivial PCA in physical domain, making these methods suffer from the common drawbacks of the trivial PCA-based ones: (1) high computational and memory demand; (2) natural structure and correlation are broken, leading to loss of potentially more compact and useful representations \cite{Lu2013}. Variable separation is indeed a promising way for reducing the computational scale of the modes and preserving the natural structure. However, existing methods can only represent random fields defined on domains which have Cartesian product decomposition format (i.e. cases where the random fields belong to \emph{structured data}), while in many cases the random fields belong to \emph{unstructured data}. Moreover, these methods generate a large number of latent factors which will bring heavy computational burden for subsequent stochastic analysis. In addition, all the existing methods in both categories requires \emph{manually} predefining tensor product basis (or collocation points) to discretize the modes and selecting quadrature points to compute the stiffness matrices. Lack of automation in these selections limits generality and adaptivity \cite{Gorodetsky2019}.

It is well known that a general non-Gaussian random field should be described with multiple statistics rather than using the covariance function only. The most widely used description which combines the covariance function and single-point marginal distribution function cannot reflect the higher-order and nonlinear correlation structure of the random field. It was pointed out in \cite{Morton2009} that joint distributions of multiple points can be non-Gaussian even if each marginal distribution is Gaussian. It is also well known that the cumulant functions with order three or higher are always zeros, making the cumulant functions be important measures of non-Gaussianity. Therefore, it will be better to develop cumulant descriptions of the latent factors based on prescribed cumulant functions. 

This paper aims at overcoming the challenges in both aspects mentioned above by proposing novel theoretical and algorithm frameworks. Motivated by the aforementioned IGA-based works, we also use isogeometric transformation to deal with complex domains. However, rather than following the trivial PCA in the physical domain, we propose a new architecture by representing the random field on the \emph{parametric domain} in \emph{tensor train}(TT) format. The core of the architecture is a newly developed \emph{rank-revealing} algorithm for separating variables of the cumulant functions. The need for predefining tensor product basis or collocation points is totally removed, which greatly enhances adaptivity. Moreover, higher-order cumulant tensors of the latent factors can be conveniently computed in a unified framework, which is beneficial for further dimension reduction.

The rest of this paper is organized in five sections. The proposed theoretical framework and the computational procedures are detailed in Section \ref{Thefram} and \ref{Compu}, respectively. In Section \ref{Compa}, we compare the proposed framework with three related methods. Particularly, the relationship with independent component analysis (ICA) will be discussed since this approach also uses higher-order cumulants for dimension reduction. In Section \ref{Cases}, three examples with increasing dimensions are employed to validate the performance of the proposed framework.

\section{Theoretical framework of the proposed method}\label{Thefram}

Let $(\Omega, \mathcal{F}, \mathbb{P})$ be a probability space where $\Omega$ is a sample space, $\mathcal{F}$ is a $\sigma$-field on $\Omega$ and $\mathbb{P}$ is a probability measure. Let $(E, \mathcal{E})$ is the measurable space on the admissible set $E$ where each $f \in E$ is defined on a bounded domain $D \subseteq \mathbb{R}^d, d=1,2,3$. A random field is defined as a measurable mapping $\omega: (\Omega, \mathcal{F}, \mathbb{P}) \to (E, \mathcal{E})$. To explicitly represent this abstract mapping, an auxiliary measurable space $(\mathbb{R}^n, \mathcal{B}^n), n\in \mathbb{N}^+$ is needed and the corresponding two mappings $\Theta: (\Omega, \mathcal{F}, \mathbb{P}) \to (\mathbb{R}^n, \mathcal{B}^n)$ and $H: (\mathbb{R}^n, \mathcal{B}^n) \to (E, \mathcal{E})$ are to be found such that $\omega=H\circ \Theta$. The task is to capture the major part of information in $ \omega $ with $ n $ as small as possible.

\subsection{Space transformation}

The key to overcome the first difficulty is to represent the physical domain with a \textit{structured} parametric domain. For a wide variety of curve-type, surface-type and solid-type domains, \textit{exact} coordinate transformations can be constructed by using NURBS-based isogeometric transformation formulated as Eq.(\ref{igfs})
\begin{subequations}\label{igfs}
    \begin{alignat}{3}
    \bm{x}(\xi) &= \sum\limits_{i=1}^{n_1} R_{i}^p (\xi) \bm{B}_i,\ &\xi\in[0,1] \quad&(\mathrm{curves})\\
    \bm{x}(\xi, \eta)& = \sum\limits_{i=1}^{n_1}\sum\limits_{j=1}^{n_2} R_{i,j}^{p,q} (\xi, \eta) \bm{B}_{i,j},\ &(\xi,\eta)\in[0,1]^2\quad &(\mathrm{surfaces})\\
    \bm{x}(\xi, \eta, \zeta)& = \sum\limits_{i=1}^{n_1}\sum\limits_{j=1}^{n_2}\sum\limits_{k=1}^{n_3} R_{i,j,k}^{p,q,r} (\xi, \eta, \zeta) \bm{B}_{i,j,k},\ &(\xi,\eta,\zeta)\in[0,1]^3\quad &(\mathrm{solids})
    \end{alignat}
\end{subequations} 
where $R_{i}^p, R_{i,j}^{p,q}$ and $R_{i,j,k}^{p,q,r}$ are NURBS basis functions defined as Eq.(\ref{nurbsR}),
\begin{subequations}\label{nurbsR}
	\begin{alignat}{1}
	 R_{i}^{p} (\xi)&= \dfrac{N_{i}^{p}(\xi)w_i}{\sum\limits_{i=1}^{n_1} N_{i}^{p}(\xi)w_i}  \label{xi1}\\
	 R_{i,j}^{p,q} (\xi,\eta)&= \dfrac{N_{i}^{p}(\xi)N_{j}^{q}(\eta)w_{i,j}}
	 {\sum\limits_{i=1}^{n_1}\sum\limits_{j=1}^{n_2} N_{i}^{p}(\xi)N_{j}^{q}(\eta)w_{i,j}} \label{xi2}\\
	 R_{i,j,k}^{p,q,r} (\xi,\eta,\zeta)&=\dfrac{N_{i}^{p}(\xi)N_{j}^{q}(\eta)N_{k}^{r}(\zeta)w_{i,j,k}}
	 {\sum\limits_{i=1}^{n_1}\sum\limits_{j=1}^{n_2}\sum\limits_{k=1}^{n_3} N_{i}^{p}(\xi)N_{j}^{q}(\eta)N_{k}^{r}(\zeta)w_{i,j,k}} \label{xi3}
	\end{alignat}
\end{subequations} 
$B_{i}, B_{i,j}$ and $B_{i,j,k}$ are control points and $w_{i}, w_{i,j}$ and $w_{i,j,k}$ are weights. Each $N_{i}^{p}$ is a $p$-degree ($(p+1)$-order) B-spline basis function. Given a knot vector $\Xi = (\xi_1,\dotsc, \xi_{n+p+1})$ (a non-decreasing set of coordinates in the parametric space [0,1]), $N_{i}^{p}$ are defined recursively by the Cox-de Boor recursion formula in Eq.(\ref{bs})
\begin{subequations}\label{bs}
	\begin{alignat}{1}
	N_{i}^0 &= \begin{cases}
	1\quad \xi\in [ \xi_i,\xi_{i+1} ) \\
	0\quad \mathrm{otherwise}
	\end{cases}\\
	N_{i}^p(\xi) &= \dfrac{\xi-\xi_{i}}{\xi_{i+p}-\xi_{i}} N_{i}^{p-1}(\xi) + \dfrac{\xi_{i+p+1}-\xi}{\xi_{i+p+1}-\xi_{i+1}} N_{i+1}^{p-1}(\xi)
	\end{alignat}
\end{subequations} 
where 0/0=0. Thus, given three ingredients: control points, weights and knot vectors, the isogeometric transformation can be constructed and formulated in a uniform format as Eq.(\ref{uniigf}).
\begin{equation}\label{uniigf}
\bm{x}(\bm{\xi}) = \sum_{\bm{I}}R_{\bm{I}}^{\bm{p}}(\bm{\xi})\bm{B_I}
\end{equation}
where $\bm{x}\in D\subseteq\mathbb{R}^d,\ \bm{\xi}\in [0,1]^{m}$ and $m\leqslant d$.

\subsection{Generalized Karhunen-Lo{\`e}ve expansion in the parametric space}\label{theocore}

After representing the physical coordinates with parametric ones, the original random field $\omega(\bm{x},\theta)$ is transformed to the parametric space. Denote $\alpha(\bm{\xi}, \theta) = \omega(\bm{x}(\bm{\xi}), \theta)$ and the corresponding covariance function $\tilde{C}(\bm{\xi},\bm{\xi}') = C(\bm{x}(\bm{\xi}),\bm{x}(\bm{\xi}'))$. We seek to represent $\alpha(\bm{\xi}, \theta)$ with a generalized K-L expansion by improving the hierarchical SVD method in \cite{Zentner2016}. Taking the case of $d$=2 as an example, the first step is to split the $\xi$-modes $\bm{f}(\xi)$ = $ (f_{i}(\xi))_{1\times n_1} $
from $ \alpha(\xi,\eta;\theta) $ by computing the eigenfunctions of the kernel in Eq.(\ref{Cf2}).
\begin{equation}\label{Cf2}
\tilde{C}_f(\xi,\xi') = \int_0^1 \tilde{C}(\xi,\eta; \xi',\eta) \mathrm{d}\eta
\end{equation}
Thus, a rank-$n_1$ decomposition is derived as Eq.(\ref{step1})
\begin{equation}\label{step1}
\alpha(\xi,\eta;\theta) \approx \bm{f}(\xi) \bm{\alpha}_{-\xi}(\eta,\theta)
\end{equation}
where  
\begin{equation}
\bm{\alpha}_{-\xi}(\eta,\theta) = \int_0^1 \bm{f}(\xi)^{\mathrm{T}}\alpha(\xi,\eta;\theta)\mathrm{d}\xi
\end{equation}
and the corresponding covariance function is a matrix-valued function in Eq.(\ref{Cg}).
\begin{equation}\label{Cg}
\tilde{C}_{-\xi}(\eta,\eta')_{n_1\times n_1} = \int_0^1 \int_0^1 \bm{f}(\xi)^{\mathrm{T}} \tilde{C}(\xi,\eta;\xi',\eta')\bm{f}(\xi')\mathrm{d}\xi\mathrm{d}\xi'
\end{equation}
Then, rather than decomposing each component of $\bm{\alpha}_{-\xi}(\eta,\theta)$ separately as in \cite{Zentner2016}, $\bm{\alpha}_{-\xi}(\eta,\theta)$ is treated as whole object by regarding the index as a new coordinate. Next, the $\eta$-modes of $ \alpha(\xi,\eta;\theta) $ (also the $\eta$-modes of $\bm{\alpha}_{-\xi}(\eta,\theta)$), denoting $\bm{g}(\eta)_{n_1\times n_2} $, are split by computing the eigenpairs of $\tilde{C}_{-\xi}(\eta,\eta')$, see Eq.(\ref{etamodes}).
\begin{equation}\label{etamodes}
\int_0^1 \tilde{C}_{-\xi}(\eta,\eta') \bm{g}_{:,j}(\eta')\mathrm{d}\eta' = \mu_j \bm{g}_{:,j}(\eta)
\end{equation}
Finally, we get the \textit{tensor train} decomposition of $\alpha(\xi,\eta;\theta)$ as Eq.(\ref{alpha2tt})
\begin{equation}\label{alpha2tt}
\alpha(\xi,\eta;\theta) \approx \bm{f}(\xi) \bm{g}(\eta) \bm{\gamma}(\theta)
\end{equation}
where 
\begin{equation}\label{gamma2}
\bm{\gamma}(\theta) = \int_0^1\int_0^1 \bm{g}(\eta)^{\mathrm{T}}\bm{f}(\xi)^{\mathrm{T}} \alpha(\xi,\eta;\theta) \mathrm{d}\xi\mathrm{d}\eta
\end{equation}
The component functions in Eq.(\ref{alpha2tt}) satisfy orthogonal conditions in Eq.(\ref{ortho1})
\begin{subequations}\label{ortho1}
	\begin{alignat}{1}
	\langle f_i, f_j \rangle &=\int_0^1 f_i(\xi) f_j(\xi) \mathrm{d}\xi = \delta_{ij} \\
	\langle \bm{g}_{:,i}, \bm{g}_{:,j} \rangle &= \int_0^1 \bm{g}_{:,i}(\eta)^{\mathrm{T}} \bm{g}_{:,j}(\eta)\mathrm{d}\eta = \delta_{ij} \\
	\langle \gamma_i, \gamma_j \rangle &=E[h_i(\theta) h_j(\theta)] = \mu_i\delta_{ij}
	\end{alignat}
\end{subequations}
 All the results can be directly extended to the case of $d$=3 as Eqs.(\ref{step2b}) to (\ref{ortho2}).
\begin{equation}\label{step2b}
\alpha(\xi,\eta, \zeta;\theta) \approx \bm{f}(\xi)_{1\times n_1} \bm{g}(\eta)_{n_1\times n_2}\bm{h}(\zeta)_{n_2\times n_3} \bm{\gamma}(\theta)_{n_3\times 1}
\end{equation}
where each column of $ \bm{f}(\xi) $ is an eigenfunction of the covariance kernel $\tilde{C}_{f}(\xi,\xi')$ in Eq.(\ref{Cf3}),
\begin{equation}\label{Cf3}
\tilde{C}_{f}(\eta,\eta')= \int_0^1 \int_0^1 \tilde{C}(\xi,\eta,\zeta;\xi',\eta,\zeta) \mathrm{d}\eta\mathrm{d}\zeta
\end{equation}
each column of $ \bm{g}(\eta) $ is an eigenfunction of $\tilde{C}_{-\xi}(\eta,\eta')$ in Eq.(\ref{Cg2})
\begin{equation}\label{Cg2}
\tilde{C}_{-\xi}(\eta,\eta')_{n_1\times n_1} = \int_0^1 \int_0^1 \bm{f}(\xi)^{\mathrm{T}} \int_0^1 \tilde{C}(\xi,\eta,\zeta;\xi',\eta',\zeta) \mathrm{d}\zeta \bm{f}(\xi')\mathrm{d}\xi\mathrm{d}\xi'
\end{equation}
and each column of $ \bm{h}(\zeta) $ is an eigenfunction of $\tilde{C}_{-\xi\eta}(\zeta,\zeta')$ in Eq.(\ref{Ch}).
\begin{equation}\label{Ch}
\begin{split}
\tilde{C}_{-\xi\eta}(\zeta,\zeta')_{n_2\times n_2} =& \int_0^1 \int_0^1 \int_0^1 \int_0^1 \bm{g}(\eta)^{\mathrm{T}} \bm{f}(\xi)^{\mathrm{T}} \tilde{C}(\xi,\eta,\zeta;\xi',\eta',\zeta')\\
&\bm{f}(\xi') \bm{g}(\eta'))\mathrm{d}\xi\mathrm{d}\eta\mathrm{d}\xi'\mathrm{d}\eta'
\end{split}
\end{equation}
\begin{equation}\label{gamma3}
\bm{\gamma}(\theta) = \int_0^1 \int_0^1\int_0^1 \bm{h}(\zeta)^{\mathrm{T}} \bm{g}(\eta)^{\mathrm{T}}\bm{f}(\xi)^{\mathrm{T}} \alpha(\xi,\eta,\zeta;\theta) \mathrm{d}\xi\mathrm{d}\eta\mathrm{d}\zeta
\end{equation}
\begin{subequations}\label{ortho2}
	\begin{alignat}{1}
	\langle f_i, f_j \rangle &=\int_0^1 f_i(\xi) f_j(\xi) \mathrm{d}\xi = \delta_{ij} \\
	\langle \bm{g}_{:,i}, \bm{g}_{:,j} \rangle &= \int_0^1 \bm{g}_{:,i}(\eta)^{\mathrm{T}} \bm{g}_{:,j}(\eta)\mathrm{d}\eta = \delta_{ij} \\
	\langle \bm{h}_{:,i}, \bm{h}_{:,j} \rangle &= \int_0^1 \bm{h}_{:,i}(\eta)^{\mathrm{T}} \bm{h}_{:,j}(\eta)\mathrm{d}\eta = \delta_{ij} \\
	\langle \gamma_i, \gamma_j \rangle &=E[\gamma_i(\theta) \gamma_j(\theta)] = \nu_i\delta_{ij}
	\end{alignat}
\end{subequations}
where $\nu_i$ is the $i$th eigenvalue of the covariance kernel $\tilde{C}_{-\xi\eta}(\zeta,\zeta')$ defined in Eq.(\ref{Ch}).

Taking $m$=2 as an example, we have the following theorem:
\begin{theorem}\label{theo1}
	The generalized K-L expansion in Eq.(\ref{alpha2tt}) is equivalent to the trivial PCA in the parametric space when $n_1,\ n_2\to\infty$.
\end{theorem}
\begin{proof}
	Denote $E_1$ = $\{F_k^{(1)}(\xi,\eta)\}_{k=1}^{\infty} $ and $\Lambda$ = $\{\lambda_k\}_{k=1}^{\infty}$ as the sets of 2D modes and eigenvalues in the trivial PCA and $E_2$ = $\{F_k^{(2)} (\xi,\eta) = \bm{f}(\xi)\bm{g}_{:,k}(\eta)\}_{k=1}^{\infty}$ and $\rm M$ = $\{\mu_k\}_{k=1}^{\infty}$ as the set of 2D modes in the generalized K-L expansion. We only need to prove that $E_1$ = $E_2$ and $\Lambda$ = $\rm M$.
	
	First, we prove $E_1 \subseteq E_2$ and $\Lambda\in \rm M$. For each $k=1,2,\dotsc$, by projecting $F_k(\xi,\eta)$ on $\bm{f}(\xi)$, 
	\begin{equation*}
	F_k(\xi,\eta) = \bm{f}(\xi)\hat{\bm{g}}_{:,k}(\eta)
	\end{equation*}
	we derive
	\begin{equation*}
	\begin{split}
	& \int_0^1 \int_0^1 \tilde{C}(\xi,\eta;\xi',\eta') \bm{f}(\xi')\hat{\bm{g}}_{:,k}(\eta') \mathrm{d}\xi'\mathrm{d}\eta' = \lambda_k \bm{f}(\xi)\hat{\bm{g}}_{:,k}(\eta)\\
	\Rightarrow &\int_0^1 \bm{f}(\xi)^{\mathrm{T}} \int_0^1 \int_0^1 \tilde{C}(\xi,\eta;\xi',\eta') \bm{f}(\xi')\hat{\bm{g}}_{:,k}(\eta') \mathrm{d}\xi'\mathrm{d}\eta' \mathrm{d}\xi = \lambda_k \int_0^1 \bm{f}(\xi)^{\mathrm{T}}\bm{f}(\xi)\mathrm{d}\xi\hat{\bm{g}}_{:,k}(\eta)\\
	\Rightarrow & \int_0^1 \tilde{C}_{-\xi}(\eta,\eta') \hat{\bm{g}}_{:,k}(\eta')\mathrm{d}\eta' = \lambda_k \hat{\bm{g}}_{:,k}(\eta)\\
	\Rightarrow & \hat{\bm{g}}_{:,k}(\eta) \mathrm{\ is\ a\ column\ of\ } \bm{g}(\eta), \lambda_k\in \rm M\\
	\Rightarrow & E_1\subseteq E_2,  \Lambda\in \rm M.
	\end{split}     
	\end{equation*}
	
	Then, we prove $E_2 \subseteq E_1$ and $\rm M\in \Lambda$. For each $k=1,2,\dotsc$, we derive
	\begin{equation*}
	\begin{split}
	&\varepsilon_k(\xi,\eta) =  \int_0^1 \int_0^1 \tilde{C}(\xi,\eta;\xi',\eta') \bm{f}(\xi')\bm{g}_{:,k}(\eta') \mathrm{d}\xi'\mathrm{d}\eta' - \mu_k \bm{f}(\xi)\bm{g}_{:,k}(\eta)\\
	\Rightarrow & \int_0^1 \bm{f}(\xi)^{\mathrm{T}} \varepsilon_k(\xi,\eta) \mathrm{d}\xi = \int_0^1 \tilde{C}_{-\xi}(\eta,\eta') \bm{g}_{:,k}(\eta')\mathrm{d}\eta' - \mu_k \bm{g}_{:,k}(\eta)=\bm{0}\\
	\Rightarrow & \varepsilon_k(\xi,\eta) \equiv 0\\
	\Rightarrow & \int_0^1 \int_0^1 \tilde{C}(\xi,\eta;\xi',\eta') \bm{f}(\xi')\bm{g}_{:,k}(\eta') \mathrm{d}\xi'\mathrm{d}\eta' = \mu_k \bm{f}(\xi)\bm{g}_{:,k}(\eta)\\
	\Rightarrow & E_2 \subseteq E_1, \rm M\in \Lambda.
	\end{split}	
	\end{equation*}
	
	Finally, we get $E_1$ = $E_2$ and $\Lambda$ = $\rm M$.
	
\end{proof}
Thus, Eq.(\ref{alpha2tt}) is mean-square convergent when $n_1,\ n_2\to\infty$. in  It is easy to verify that this theorem is valid for other values of $m$.

\subsection{Representing higher-order cumulants of latent factors}

It is well known that the probabilistic structure of a random field is uniquely defined by its family of finite-dimensional marginal distributions. For arbitrary $n$ points $\bm{\Xi}=(\bm{\xi}_1,\dotsc,\bm{\xi}_n)$ in parametric space and arbitrary $\bm{u}=(u_1,\dotsc,u_n)\in\mathbb{R}^n$, according to the A-type Gram-Charlier series, the corresponding marginal probability density function (PDF) can be represented as Eq.(\ref{AGC}) 
\begin{equation}\label{AGC}
\begin{split}
f_{\bm{\Xi}}(\bm{u};\bm{\kappa}) =& f_G(\bm{u}) \left[1 + \dfrac{1}{3!}\sum\limits_{i,j,k} \kappa^{i,j,k}h_{ijk}(\bm{u})+ \dfrac{1}{4!}\kappa^{i,j,k,l}h_{ijkl}(\bm{u}) +\dfrac{1}{5!}\kappa^{i,j,k,l,m}h_{ijklm}(\bm{u})+ \right.\\
&\left.\dfrac{1}{6!} \left(\kappa^{i,j,k,l,m,n}+10\kappa^{i,j,k}\kappa^{l,m,n}\right) h_{ijklmn}(\bm{u})+ \cdots \right]
\end{split}
\end{equation}
where $f_G$ is the PDF of the Gaussian distribution with mean zero and covariance function $\tilde{C}(\bm{\xi},\bm{\xi}')$, $h$ with subscripts are Hermite polynomials and $\bm{\kappa}$ is the family of cumulant functions consists of elements $\kappa^{i,j,k}$ = $\tilde{C}_3(\bm{\xi}_i,\bm{\xi}_j,\bm{\xi}_k) $ (the third-order cumulant function) and so on. Meanwhile, $\bm{\kappa}$ is defined by the family of marginal PDFs by definition. Thus, $\bm{\kappa}$ is equivalent to the family of marginal PDFs by definition in terms of describing the probabilistic structure of a random field. In practical problems, only finite orders of cumulant functions are available. To absorb the information of high-order ($\geqslant$3) cumulant functions into the representation of $\alpha(\bm{\xi};\theta)$, taking the case of $d$=2 as an example, since the reconstructed random field has the form in Eq.(\ref{alpha2tt}),  the relationship between the third-order cumulant functions and the third-order cumulants of latent factors is expressed as Eq.(\ref{C3}).
\begin{equation}\label{C3}
\begin{split}
\breve{C}_3(i,j,k) =& \int_{[0,1]^6} \tilde{C}_3(\xi,\eta;\xi',\eta';\xi'',\eta'') (\bm{g}(\eta)^{\mathrm{T}}_{:,i}\bm{f}(\xi)^{\mathrm{T}}) (\bm{g}(\eta')^{\mathrm{T}}_{:,j}\bm{f}(\xi')^{\mathrm{T}})\\
& (\bm{g}(\eta'')^{\mathrm{T}}_{:,k}\bm{f}(\xi'')^{\mathrm{T}})\mathrm{d}\xi\mathrm{d}\eta
\mathrm{d}\xi'\mathrm{d}\eta'\mathrm{d}\xi''\mathrm{d}\eta''
\end{split}
\end{equation}
This relationship can be directly extended to the cases of other values of $d$ and cumulant orders. Unfortunately, $ \breve{C}_3(i,j,k) $ is generally not orthogonal decomposable \cite{Robeva2016}, hence, we have to resort to HOSVD to further reduce the dimensionality of the latent factors. Finally, by combing the coordinate transformation in Eq.(\ref{uniigf}) with Eqs.(\ref{alpha2tt}) (or (\ref{step2b})) and (\ref{C3}), we get a parametric representation of the original random field $\omega(\bm{x};\theta)$ in the sense of given cumulants.

\section{Computational procedure}\label{Compu}

In this section, we seek an algorithm framework which is general enough to overcome the difficulties in both aspects mentioned in Section \ref{Intro}.

\subsection{Space transformation}

Given three ingredients: control points, weights and knot vectors, to avoid redundant computations, geometry of the physical domain defined in Eq.(\ref{uniigf}) is represented with the fast isogeometric transformation algorithm \cite{Piegl2012}.

\subsection{Generalized Karhunen-Lo{\`e}ve expansion in the parametric space}

After representing the physical domain $D$ with a structured parametric domain $[0,1]^m$, the core of the new algorithm framework is representing the modes of $\alpha(\bm{\xi};\theta)$. This task is accomplished with the following two steps.

\subsubsection{Tensor train decomposition of cumulant functions}

 According to the theoretical framework in the previous section, it will be very beneficial for fast computations if the cumulant functions have separable forms. Borrowing the idea in \cite{Savostyanov2014}, we propose a tensor train decomposition algorithm to obtain a low-rank separable approximation to a $K$th ($K\geqslant2$) cumulant function $\tilde{C}_K(\bm{\xi}_1,\dotsc,\bm{\xi}_K)$.
 
 For the sake of simplicity, $\tilde{C}_K(\bm{\xi}_1,\dotsc,\bm{\xi}_K)$ is redefined as an auxiliary function $G: [0,1]^{a}\to\mathbb{R}\ (a=mK)$ where $\bm{u}=(u_1,\dotsc,u_a) \in [0,1]^a$ is a permutation of the original coordinates. The idea is reconstruct $G$ only by using some of its fibers. More precisely, we seek a separable representation of $G$ as Eq.(\ref{sepG})
 \begin{equation}\label{sepG}
 \begin{split}
 G(u_1,\dotsc,u_a) \approx & G(u_1, \mathcal{I}^{>1})[G(\mathcal{I}^{\leqslant 1}, \mathcal{I}^{>1})]^{-1} G(\mathcal{I}^{\leqslant 1}, u_2, \mathcal{I}^{>2})\\
 &[G(\mathcal{I}^{\leqslant 2}, \mathcal{I}^{>2})]^{-1} \cdots G(\mathcal{I}^{\leqslant a-1}, u_a)
 \end{split} 
 \end{equation}
 where $\mathcal{I}^{\leqslant k}$ and $\mathcal{I}^{>k}$ are the $k$th pair of interpolation sets (both have cardinality $r_k$ which is the $k$th rank). The optimal choice of the $k$th pair of interpolation sets is the solution to Eq.(\ref{maxvol}).
 \begin{equation}\label{maxvol}
 \{\mathcal{I}^{\leqslant k*}, \mathcal{I}^{>k*}\} =\arg\max_{\substack{\mathcal{I}^{\leqslant k}, \mathcal{I}^{>k}}} |\det G(\mathcal{I}^{\leqslant k}, \mathcal{I}^{>k})|
 \end{equation}
However, the search for the maximum-volume submatrix is an NP (non-deterministic polynomial)-hard problem. Hence, we propose a \emph{heuristic} algorithm to find a quasi-optimal choice of $\{\mathcal{I}^{\leqslant k}, \mathcal{I}^{>k}\}_{k=1}^{a-1}$ and derive an adaptive tensor train decomposition of a cumulant function, see Algorithm \ref{ctt}.
\begin{algorithm}
	\caption{Adaptive tensor train decomposition of a $K$th cumulant function on a $m$-dimensional parametric space}\label{ctt}
	\begin{algorithmic}[1]
	\Require  Auxiliary function $G: [0,1]^{a}\to\mathbb{R}\ (a=mK)$; maximum number of sweeps $maxswp$; stopping tolerance $tol$
	\Ensure Tensor train decomposition of $G$, denoting $G_{\mathrm{TT}}$; interpolation sets $\mathcal{I}$; a matrix $errdm$ containing iteration errors of each dimension during each iteration	
	
	\State $\bm{u}^{(0)} \gets \arg\max_{\substack{\bm{u}}} |G(\bm{u})|$ from $m_0$ quasi-random samples;\Comment{Find the initial pivot}	
	\ForAll{$k = 1: a-1$}
	\State{$\{\mathcal{I}^{\leqslant k}, \mathcal{I}^{>k}\} \gets \{\{\bm{u}^{(0)\leqslant k} \}, \{\bm{u}^{(0)> k} \}\}$;}\Comment{Initialize the interpolation sets}
	\EndFor
	\State $\mathcal{I} \gets \{\mathcal{I}^{\leqslant k}, \mathcal{I}^{>k}\}_{k=1}^{a-1}$;	
	\ForAll{$k=1:a-1$}
	\State{$errd(k)=1$;}\Comment{Initialize the iteration error in each dimension}
	\EndFor	
	 \State $S\gets$ 0;\Comment{Initialize the sweep number}
	 \While{$S<maxswp$}
	 \State $S \gets S+1$;\Comment{Begin left-to-right sweep}	
	\For{$k=1:a-1$}
		\If{$ errd(k)<tol $}
			\State {continue}
		\Else		
		   \State $ [ \mathcal{I} ,  errmax ]\gets $\Call{ISE}{ $ k $ ,$ \mathcal{I} $,$ G $, $ a-1 $ } \Comment{use Algorithm \ref{ISE} to expand interpolation sets }
		   \State $errd(k)\gets errmax$;\Comment{update iteration errors}
	   \EndIf
	\EndFor    
	\If{$S=1$}\Comment{Record iteration errors of all dimensions during each iteration}
	\State{$errdm\gets errd$}
	\Else
	\State {$errdm\gets [errdm,errd]$}	
	\EndIf
	\If{$\max(errd)<tol$ }
	\State {break}
	\EndIf
	
	\State $S\gets S+1$;\Comment{Begin right-to-left sweep}
   \For{$k=a:-1:2$}
   	\If{$ errd(k-1)<tol $}
   	\State{continue}
   	\Else   		
   		\State $[\mathcal{I}, errmax]\gets$ \Call{ISE}{$ k -1$,$\mathcal{I}$,$G$, $ a-1 $}\Comment{use Algorithm \ref{ISE} to expand interpolation sets }
   		\State $errd(k-1)\gets errmax$;\Comment{update iteration errors}
   \EndIf
   \EndFor  
   \State $errdm\gets [errdm,errd]$; \Comment{Record iteration errors of all dimensions during each iteration}
   \algstore{bkbreak}
   \end{algorithmic}
\end{algorithm}
\begin{algorithm}
	\begin{algorithmic}
		\algrestore{bkbreak}
   \If{$\max(errd)<tol$}
   \State {break}		
   \EndIf
   \EndWhile
\ForAll{$k=1:a$}\Comment{Output tensor train decomposition by adaptive reconstruction of each fiber}
\If{$k=1$}
\State $\hat{\bm{G}}(u_1)_{1\times r_1}\gets$   \Call{chebmatrix} {$G( u_1 , \mathcal{I}^{>1} )$}  $[G(\mathcal{I}^{\leqslant 1}, \mathcal{I}^{>1})]^{-1} $; \label{chb1}
\ElsIf{$k<d$}
\State $\hat{\bm{G}}(u_k)_{r_{k-1}\times r_k}\gets$  \Call{chebmatrix}{ $ G(\mathcal{I}^{\leqslant k-1} , u_k, \mathcal{I}^{>k})$ } $ [G(\mathcal{I}^{\leqslant k}, \mathcal{I}^{>k})]^{-1} $; \label{chb2}
\Else
\State $\hat{\bm{G}}(u_a)_{r_{a-1}\times 1} \gets$ \Call{chebmatrix} {$ G(\mathcal{I}^{\leqslant a-1}, u_a) $};\label{chb3}
\EndIf
\EndFor
\State\begin{equation}
G_{\mathrm{TT}}(\bm{u})= \hat{\bm{G}}_1(u_1) \hat{\bm{G}}_2(u_2) \cdots \hat{\bm{G}}_a(u_a)
\end{equation}
\end{algorithmic}
\end{algorithm}
After initialization, the interpolation sets for each two neighboring dimensions are progressively enriched by the left-to-right and right-to-left sweeps until the maximum approximation error is smaller than the prescribed threshold. The interpolation sets enriched in this way are two-side nested, i.e. $\forall k\in\{1,\dotsc,a-1\}$, $\mathcal{I}^{\leqslant k+1}\subseteq \mathcal{I}^{\leqslant k}\times [0,1]$, $\mathcal{I}^{> k}\subseteq [0,1]\times\mathcal{I}^{> k+1}$. Then, each fiber in Eq.(\ref{sepG}) is adaptively represented in the \emph{chebfun} format \cite{Driscoll2014}. Next, fibers of each two neighboring dimensions form a \emph{chebmatrix} as in lines \ref{chb1}, \ref{chb2}, \ref{chb3}.
\begin{algorithm}
	\caption{ $ ISE $: Interpolation set expansion algorithm}\label{ISE}
	\begin{algorithmic}[1]
	\Require {current dimension $k$; current interpolation sets $\mathcal{I}$; auxiliary function $ G $; number of variables of $G $, denoting $a$; $m_k$ for solving Eq.(\ref{newpiv}); iteration error tolerance $ tol $ }
	\Ensure {Expanded interpolation sets $\mathcal{I}$; iteration error $errmax$}
	\If{$k$=1}
	\State Join the supports of the matrix-valued function $G(u_k, u_{k+1}, \mathcal{I}^{>k+1})$ in the second dimension and define a new function $G_{\mathrm{Jst}}: [0,1]\times [0,r_{k+1}]\to\mathbb{R}$;
	\State Redefine the matrix-valued function $G(u_k,  \mathcal{I}^{>k})$ as a new function $G_{\mathrm{Js}}: [0,1]\to\mathbb{R}^{1\times r_k}$;
	\State Join the supports of each row of the matrix-valued function $G(\mathcal{I}^{\leqslant k}, u_{k+1}, \mathcal{I}^{>k+1})$ and define a new function $G_{\mathrm{Jt}}: [0,r_{k+1}]\to\mathbb{R}^{r_k\times 1}$;		
	\ElsIf{$k<d-1$}
	\State Join the supports of the matrix-valued function $G( \mathcal{I}^{\leqslant k-1}, u_k, u_{k+1}, \mathcal{I}^{>k+1})$ and define a new function $G_{\mathrm{Jst}}: [0,r_{k-1}]\times [0,r_{k+1}]\to\mathbb{R}$;
	\State Join the supports of each column of the matrix-valued function $G(\mathcal{I}^{\leqslant k-1}, u_k,  \mathcal{I}^{>k})$ and define a new function $G_{\mathrm{Js}}: [0,r_{k-1}]\to\mathbb{R}^{1\times r_k}$;
	\State Join the supports of each row of the matrix-valued function $G(\mathcal{I}^{\leqslant k}, u_{k+1}, \mathcal{I}^{>k+1})$ and define a new function $G_{\mathrm{Jt}}: [0,r_{k+1}]\to\mathbb{R}^{r_k\times 1}$;
	\Else
	\State Join the supports of the matrix-valued function $G(\mathcal{I}^{\leqslant k-1}, u_k, u_{k+1})$ in the first dimension and define a new function $G_{\mathrm{Jst}}: [0,r_{k-1}]\times [0,1]\to\mathbb{R}$;
	\State Join the supports of each column of the matrix-valued function $G(\mathcal{I}^{\leqslant k-1}, u_k,  \mathcal{I}^{>k})$ and define a new function $G_{\mathrm{Js}}: [0,r_{k-1}]\to\mathbb{R}^{1\times r_k}$;
	\State Redefine the matrix-valued function $G(\mathcal{I}^{\leqslant k}, u_{k+1})$ as a new function $G_{\mathrm{Jt}}: [0,1]\to\mathbb{R}^{r_k\times 1}$;
	\EndIf
\State Compute $[G(\mathcal{I}^{\leqslant k}, \mathcal{I}^{>k})]^{-1}$;
\State Solve Eq.(\ref{newpiv})
\begin{equation}\label{newpiv}
(s^*, t^*) = \arg\max_{\substack{s,t}} |G_{\mathrm{Jst}}(s,t) - G_{\mathrm{Js}}(s) [G(\mathcal{I}^{\leqslant k}, \mathcal{I}^{>k})]^{-1} G_{\mathrm{Jt}}(t)|
\end{equation}
with  $m_k$ quasi-random samples;
\State $errmax \gets |G_{\mathrm{Jst}}(s^*,t^*) - G_{\mathrm{Js}}(s^*) [G(\mathcal{I}^{\leqslant k}, \mathcal{I}^{>k})]^{-1} G_{\mathrm{Jt}}(t^*)|$;\Comment{Iteration error}
\If{$ errmax<tol $}
\State break
\EndIf
\State $i^*\gets [s^*]$, $j^*\gets  [t^*]$\Comment{$[x]$ denotes the largest integer no larger than $x$}
\State $u^{(S)}$ is obtained with $u^{(S)\leqslant k} \gets [\mathcal{I}^{\leqslant k-1}(i^*,:), s^*-i^*]$ and $u^{(S)> k} \gets [t^*-j^*, \mathcal{I}^{>k+1}(j^*,:)]$;\Comment{The new pivot}	
\State $ \mathcal{I}^{\leqslant k} \gets  \mathcal{I}^{\leqslant k} \cup \{u^{(S)\leqslant k}\} $, $ \mathcal{I}^{>k}  \gets  \mathcal{I}^{>k} \cup \{u^{(S)> k}\}$;\Comment{The expanded interpolation set of dimension $k$}
\end{algorithmic}
\end{algorithm}

For a special case where $K=2$ and $m=1$, Algorithm \ref{ctt} degenerates into Algorithm \ref{ctt2}.

\begin{algorithm}
	\caption{Adaptive tensor train decomposition of a bi-variate covariance function}\label{ctt2}
\begin{algorithmic}[1]
	\Require  Auxiliary function $G: [0,1]^{2}\to\mathbb{R}$; maximum number of sweeps $maxswp$; stopping tolerance $tol$
	\Ensure Tensor train decomposition of $G$, denoting $G_{\mathrm{TT}}$; interpolation sets $\mathcal{I}$; a matrix $errdm$ containing iteration errors during each iteration	
	
	\State $\bm{u}^{(0)} \gets \arg\max_{\substack{\bm{u}}} |G(\bm{u})|$ from $m_0$ quasi-random samples;\Comment{Find the initial pivot}	
	\State{$\{\mathcal{I}^{\leqslant 1}, \mathcal{I}^{>1}\} \gets \{\{\bm{u}^{(0)\leqslant 1} \}, \{\bm{u}^{(0)> 1} \}\}$;}\Comment{Initialize the interpolation sets}
	\State $\mathcal{I} \gets \{\mathcal{I}^{\leqslant 1}, \mathcal{I}^{>1}\}$;	
	\State{$errmax=1$;}\Comment{Initialize the iteration error}
	\For{$S=1:maxswp$}
	\State Compute $[G(\mathcal{I}^{\leqslant 1}, \mathcal{I}^{>1})]^{-1}$;
	\State\label{m1} Solve Eq.(\ref{newpiv_})
	\begin{equation}\label{newpiv_}
	(s^*, t^*) = \arg\max_{\substack{s,t}} |G(s,t) - G(s,  \mathcal{I}^{>1}) [G(\mathcal{I}^{\leqslant 1}, \mathcal{I}^{>1})]^{-1} G(\mathcal{I}^{\leqslant k}, t)|
	\end{equation}
	with  $m_1$ quasi-random samples;
	\State $errmax \gets |G(s^*,t^*) - G(s^*,  \mathcal{I}^{>1}) [G(\mathcal{I}^{\leqslant 1}, \mathcal{I}^{>1})]^{-1} G(\mathcal{I}^{\leqslant k}, t^*)|$;\Comment{Iteration error}	
	\If{$S=1$}\Comment{Record iteration errors of all dimensions during each iteration}
	\State{$errdm\gets errmax$}
	\Else
	\State {$errdm\gets [errdm,errmax]$}	
	\EndIf
	\If{$errmax<tol$}
	\State {break}
	\EndIf
	\State $u^{(S)}$ is obtained with $u^{(S)\leqslant 1} \gets  s^*$ and $u^{(S)> k} \gets t^*$;\Comment{The new pivot}	
	\State $ \mathcal{I}^{\leqslant 1} \gets  \mathcal{I}^{\leqslant 1} \cup \{u^{(S)\leqslant 1}\} $, $ \mathcal{I}^{>1}  \gets  \mathcal{I}^{>1} \cup \{u^{(S)> 1}\}$;\Comment{The expanded interpolation set}	
	\EndFor
	
	\State $\hat{\bm{G}}(u_1)_{1\times r_1}\gets$   \Call{chebmatrix} {$G( u_1 , \mathcal{I}^{>1} )$}  $[G(\mathcal{I}^{\leqslant 1}, \mathcal{I}^{>1})]^{-1} $; \Comment{Output tensor train decomposition by adaptive reconstruction of each fiber}
	\State $\hat{\bm{G}}(u_2)_{r_{1}\times 1} \gets$ \Call{chebmatrix} {$ G(\mathcal{I}^{\leqslant 1}, u_2) $};
	\State\begin{equation}
	G_{\mathrm{TT}}(\bm{u})= \hat{\bm{G}}_1(u_1) \hat{\bm{G}}_2(u_2)
	\end{equation}
\end{algorithmic}	
\end{algorithm}

\subsubsection{Computation of modes in each direction using reconstructed covariance functions}

 Computation of modes in each direction is performed with $K=2$. When $m=1$, decomposition of the $\tilde{C}(\xi,\xi')$ is trivial by using Algorithm \ref{ctt2}, and we obtain a separable representation of $\tilde{C}(\xi,\xi')$ as Eq.(\ref{C2}).
 \begin{equation}\label{C2}
 \tilde{C}_{\mathrm{TT}}(\xi,\xi') = \bm{M}_1(\xi)_{1\times r_1}\bm{M}_2(\xi')_{r_1\times 1}
 \end{equation}
Before computing the modes, accuracy of the TT decomposition should be checked by computing the empirical global relative error in Eq.(\ref{epg})
\begin{equation}\label{epg}
\varepsilon_{\mathrm{g}} = \dfrac{\sqrt{\dfrac{1}{N}\sum\limits_{i=1}^{N}(\tilde{C}(\xi_i, \xi'_i) - \tilde{C}_{\mathrm{TT}}(\xi_i, \xi'_i))^2 }}{\sqrt{\dfrac{1}{N}\sum\limits_{i=1}^{N}\tilde{C}(\xi_i, \xi'_i)^2 }} 
\end{equation}
where $\{(\xi_i,\xi'_i)\}_{i=1}^N$ is the set of test samples. If $ \varepsilon_{\mathrm{g}} $ is larger than a prescribed tolerance $ tol_{\mathrm{g}} $, then go to the next step; else, adjust values of $ m_1 $ or $ tol $ in Algorithm \ref{ctt2} and re-compute the TT decomposition. Next,taking the multiplicative form in Eq.(\ref{C2}) into account, by following the idea of \cite{Townsend2013}, we propose a specific algorithm for computing the $\xi$-modes $\bm{f}(\xi)_{1\times n_1}$  and the eigenvalues $\{\lambda_i\}_{i=1}^{n_1}$ see Algorithm \ref{svdre1}. 
 \begin{algorithm}
 	\caption{A specific algorithm for computing the eigenpairs of a bi-variate covariance kernel with multiplicative form}\label{svdre1}
 	\begin{algorithmic}[1]
 		\Require Two components of the covariance kernel $\bm{M}_1(\xi)_{1\times n}$ and $\bm{M}_2(\xi')_{n\times 1}$ in $chebmatrix$ format;
 		\Ensure Eigenvalues $\{\lambda_k\}_{k=1}^{n}$ and modes $\bm{f}(\xi)_{1\times n}$ in $chebmatrix$ format
 		
 		\State $[\bm{Q}_L(\xi)_{1\times n}, \bm{R}_L]\gets \Call{QR}{\bm{M}_1(\xi)}$;
 		\State $[\bm{Q}_R(\xi')_{1\times n}, \bm{R}_R]\gets \Call{QR}{\bm{M}_2(\xi')}$;
 		\State $[\bm{U},\bm{S},\bm{V}]\gets \Call{svd}{\bm{R}_L \bm{DR}_R^{\mathrm{T}}}$;
 		\State $\{\lambda_k\}_{k=1}^{n}\gets \mathrm{diag}(\bm{S}),\ \bm{f}(\xi)_{1\times n}\gets \bm{Q}_L(\xi)\bm{U}$;
 		
 	\end{algorithmic}
 \end{algorithm}
 For the QR decomposition in Algorithm \ref{svdre1}, the first choice is the Householder triangularization \cite{Trefethen2009}(the default choice in the subsequent numerical experiments) due to its good numerical stability. Unfortunately, this method is too time-consuming for large-scale chebmatrices since the computational time of both the plus and inner-product operations are non-negligible. Hence, a second choice is the Cholesky decomposition of $\bm{M}^{\mathrm{T}}\bm{M}$. This method produces almost the same results as the first one with much higher efficiency, and still have good numerical stability for sufficiently large eigenvalues. When $m=2$, first, the auxiliary function $G$ is defined as: $G(u_1,\dotsc,u_4) $ = $ G(\xi,\eta;\eta',\xi')$ = $ C(\bm{x}(\xi,\eta);\bm{x}_{\mathrm{r}}(\eta',\xi'))$ where $\bm{x}_{\mathrm{r}}(\eta',\xi')$ = $\bm{x}(\xi',\eta')$, i.e. $ G(\xi,\eta;\eta',\xi')$ = $\tilde{C}(\xi,\eta;\xi',\eta')$. Then, by applying Algorithm \ref{ctt} on $G$, the tensor train decomposition of $ \tilde{C}(\xi,\eta;\xi',\eta')$ is obtained as Eq.(\ref{C2r2}).
 \begin{equation}\label{C2r2}
 \begin{split}
 \tilde{C}(\xi,\eta;\xi',\eta') &\approx \tilde{C}_{\mathrm{TT}}(\xi,\eta;\xi',\eta') \\
 &= \hat{\bm{G}}_1(\xi)_{1\times r_1} \hat{\bm{G}}_2(\eta)_{r_1\times r_2} \hat{\bm{G}}_3(\eta')_{r_2\times r_3} \hat{\bm{G}}_4(\xi')_{r_3\times 1}
 \end{split}
 \end{equation}
After passing the global random test, covariance kernel in Eq.(\ref{Cf2}) is obtained as Eq.(\ref{partr2}).
\begin{equation}\label{partr2}
\begin{split}
\int_{0}^1 \tilde{C}(\xi,\eta;\xi',\eta) \mathrm{d}\eta &\approx \left(\hat{\bm{G}}_1(\xi) \int_0^1 \hat{\bm{G}}_2(\eta) \hat{\bm{G}}_3(\eta)\mathrm{d}\eta\right) \hat{\bm{G}}_4(\xi')\\
&=\bm{M}_1(\xi)_{1\times r_3}\bm{M}_2(\xi')_{r_3\times 1}
\end{split}
\end{equation}
 The $\xi$-modes $\bm{f}(\xi)_{1\times n_1}$ are computed by using Algorithm \ref{svdre1}. Next, the covariance kernel in Eq.(\ref{Cg}) is obtained as Eq.(\ref{Cga}).
\begin{equation}\label{Cga}
\begin{split}
\tilde{C}_{-\xi}(\eta,\eta')_{n_1\times n_1} &\approx \int_0^1 \bm{f}(\xi)^{\mathrm{T}}\hat{\bm{G}}_1(\xi) \mathrm{d}\xi  \hat{\bm{G}}_2(\eta) \hat{\bm{G}}_3(\eta') \int_0^1 \hat{\bm{G}}_4(\xi')\bm{f}(\xi') \mathrm{d}\xi'\\
&= \bm{M}_3(\eta)_{n_1\times r_2}\bm{M}_4(\eta')_{r_2\times n_1}
\end{split}
\end{equation}
To compute the $\eta$-modes $ \bm{g}(\eta) $ in Eq.(\ref{alpha2tt}), we join the supports of $ \tilde{C}_{-\xi}(\eta,\eta') $ in both dimensions, which is equivalent to joining each column of $ \bm{M}_3(\eta) $ and each row of $ \bm{M}_4(\eta') $, respectively, deriving Eq.(\ref{gJ}).
\begin{equation}\label{gJ}
\tilde{C}_{-\xi\mathrm{J}}(\eta_{\mathrm{J}},\eta_{\mathrm{J}}') \approx \bm{M}_{3\mathrm{J}}(\eta_{\mathrm{J}})_{1\times r_2}
\bm{M}_{4\mathrm{J}}(\eta_{\mathrm{J}}')_{r_2\times 1}
\end{equation}
Then, by applying Algorithm \ref{svdre1} on $ \bm{M}_{3\mathrm{J}}(\eta_{\mathrm{J}}) $ and $ \bm{M}_{4\mathrm{J}}(\eta_{\mathrm{J}}') $, we obtain a joint form of $ \bm{g}(\eta) $, denoting $ \bm{g}_{\mathrm{J}}(\eta)_{1\times n_2} $, and the eigenvalues $\{\mu_i \}_{i=1}^{n_2}$. Next, $ \bm{g}(\eta)_{n_1\times n_2} $ is obtained by disjoining the support of each column of $\bm{g}_{\mathrm{J}}(\eta)$ with $n_1$ equal intervals.

The above procedure can be straightforwardly extended to the case of $m=3$. By letting $ G(\xi,\eta,\zeta;\eta',\xi',\zeta')$ = $\tilde{C}(\xi,\eta,\zeta;\xi',\eta',\zeta')$, first, we get Eq.(\ref{C2r3}).
 \begin{equation}\label{C2r3}
\tilde{C}(\xi,\eta,\zeta;\xi',\eta',\zeta') \approx \hat{\bm{G}}_1(\xi)_{1\times r_1} \hat{\bm{G}}_2(\eta)_{r_1\times r_2}
\hat{\bm{G}}_3(\zeta)_{r_2\times r_3}
\hat{\bm{G}}_4(\zeta')_{r_3\times r_4}
 \hat{\bm{G}}_5(\eta')_{r_4\times r_5} \hat{\bm{G}}_6(\xi')_{r_5\times 1}
\end{equation}
Then, the covariance kernel in Eq.(\ref{Cf3}) is derived as Eq.(\ref{Cf3a}),
\begin{equation}\label{Cf3a}
\begin{split}
\tilde{C}_f(\xi,\xi') &= \int_{0}^1\int_{0}^1 \tilde{C}(\xi,\eta,\zeta;\xi',\eta,\zeta) \mathrm{d}\eta\mathrm{d}\zeta \\
&\approx \left(\hat{\bm{G}}_1(\xi) \int_0^1 \hat{\bm{G}}_2(\eta) \int_0^1 \hat{\bm{G}}_3(\zeta) \hat{\bm{G}}_4(\zeta)\mathrm{d}\zeta \hat{\bm{G}}_5(\eta)\mathrm{d}\eta\right) \hat{\bm{G}}_6(\xi')\\
&=\bm{M}_1(\xi)_{1\times r_5}\bm{M}_2(\xi')_{r_5\times 1}
\end{split}
\end{equation}
and the $\xi$-modes $\bm{f}(\xi)_{1\times n_1}$ are computed by applying Algorithm \ref{svdre1} on the $ \bm{M}_1(\xi) $ and $ \bm{M}_2(\xi') $ above. Next, the covariance kernel in Eq.(\ref{Cg2}) is derived as Eq.(\ref{Cg2a})
\begin{equation}\label{Cg2a}
\begin{split}
\tilde{C}_{-\xi}(\eta,\eta')_{n_1\times n_1} =& \int_0^1 \int_0^1 \bm{f}(\xi)^{\mathrm{T}} \int_0^1 \tilde{C}(\xi,\eta,\zeta;\xi',\eta',\zeta) \mathrm{d}\zeta \bm{f}(\xi')\mathrm{d}\xi\mathrm{d}\xi'\\
\approx & \left( \int_0^1 \bm{f}(\xi)^{\mathrm{T}}\hat{\bm{G}}_1(\xi)\mathrm{d}\xi \hat{\bm{G}}_2(\eta) \int_0^1 \hat{\bm{G}}_3(\zeta) \hat{\bm{G}}_4(\zeta)\mathrm{d}\zeta \right)\\ &\ \left(\hat{\bm{G}}_5(\eta') \int_0^1\hat{\bm{G}}_6(\xi')\bm{f}(\xi')\mathrm{d}\xi'\right)\\
=&\ \bm{M}_3(\eta)_{n_1\times r_4}\bm{M}_4(\eta')_{r_4\times n_1}
\end{split}
\end{equation}
By joining supports, applying Algorithm \ref{svdre1} and disjoining supports sequentially, the $\eta$-modes $\bm{g}(\eta)$ are obtained. Next, the covariance kernel in Eq.(\ref{Ch}) is derived as Eq.(\ref{Cha}),
\begin{equation}\label{Cha}
\begin{split}
\tilde{C}_{-\xi\eta}(\zeta,\zeta')_{n_2\times n_2} =&\ \int_0^1 \int_0^1 \int_0^1 \int_0^1 \bm{g}(\eta)^{\mathrm{T}} \bm{f}(\xi)^{\mathrm{T}} \tilde{C}(\xi,\eta,\zeta;\xi',\eta',\zeta')\\
&\ \bm{f}(\xi') \bm{g}(\eta'))\mathrm{d}\xi\mathrm{d}\eta\mathrm{d}\xi'\mathrm{d}\eta' \\
\approx&\ \left( \int_0^1 \bm{g}(\eta)^{\mathrm{T}} \int_0^1 \bm{f}(\xi)^{\mathrm{T}} \hat{\bm{G}}_1(\xi)\mathrm{d}\xi \hat{\bm{G}}_2(\eta)\mathrm{d}\eta \hat{\bm{G}}_3(\zeta)\right)\\
&\ \left( \hat{\bm{G}}_4(\zeta) \int_0^1 \hat{\bm{G}}_5(\eta') \int_0^1 \hat{\bm{G}}_6(\xi')\bm{f}(\xi')\mathrm{d}\xi' \bm{g}(\eta')\mathrm{d}\eta' \right)\\
=&\ \bm{M}_5(\zeta)_{n_2\times r_3}\bm{M}_6(\zeta')_{r_3\times n_2}
\end{split}
\end{equation}
By joining supports, applying Algorithm \ref{svdre1} and disjoining supports sequentially, the $\zeta$-modes $\bm{h}(\zeta)$ and the eigenvalues $\{\nu_i\}_{i=1}^{n_3}$ are obtained.

\subsection{Representing higher-order cumulants of latent factors}

Having obtained the modes in all directions and the covariance matrix of the latent factors, the final task is to efficiently represent higher-order cumulant tensors of the latent factors. Taking $K$ = 3 as an example, by recalling Eqs.(\ref{gamma2}) and (\ref{gamma3}), a single latent factor is expressed as Eq.(\ref{gammaeq}) or uniformly written as Eq.(\ref{gammauni}).
\begin{subequations}\label{gammaeq}
\begin{alignat}{2}
\gamma_i(\theta) &= \int_0^1 f_i(\xi)\alpha(\xi,\theta)\mathrm{d}\xi,\quad  &(m=1)\\
\gamma_i(\theta) &= \int_0^1 \int_0^1 (\bm{g}(\eta))^{\mathrm{T}}_{i,:}\bm{f}(\xi)^{\mathrm{T}}\alpha(\xi,\eta; \theta)\mathrm{d}\xi\mathrm{d}\eta,\quad  &(m=2)\\
\gamma_i(\theta) &= \int_0^1 \int_0^1 \int_0^1 (\bm{h}(\zeta))^{\mathrm{T}}_{i,:}\bm{g}(\eta)^{\mathrm{T}}\bm{f}(\xi)^{\mathrm{T}}\alpha(\xi,\eta,\zeta; \theta)\mathrm{d}\xi\mathrm{d}\eta\mathrm{d}\zeta,\quad  &(m=3)
\end{alignat}
\end{subequations}
\begin{equation}\label{gammauni}
\gamma_i(\theta) = \int_{[0,1]^m} F_i(\bm{\xi})\alpha(\bm{\xi};\theta)\mathrm{d}\bm{\xi},\quad (m=1,2,3)
\end{equation}
By using Eq.(\ref{gammauni}), the third cumulant tensor of the latent factors is derived as Eq.(\ref{C3ijk}).
\begin{equation}\label{C3ijk}
\breve{C}_3(i,j,k) = \int_{[0,1]^{3m}} F_i(\bm{\xi})F_j(\bm{\xi}')F_k(\bm{\xi}'')\tilde{C}_3(\bm{\xi},\bm{\xi}',\bm{\xi}'')\mathrm{d}\bm{\xi}\mathrm{d}\bm{\xi}'\mathrm{d}\bm{\xi}'',\quad (m=1,2,3)
\end{equation}
By letting the auxiliary function $G\equiv\tilde{C}_3$, a tensor train decomposition of $\tilde{C}_3$ (denoted as $\tilde{C}_{3,\mathrm{TT}}$) is obtained with Algorithm \ref{ctt}. Since this approximation is multiplicatively separable, taking $m$=3 as an example, for each index $I\ (I=i,j,k)$, the corresponding part in the decomposition has the form $\hat{\bm{G}}_{3I-2}(\xi)\hat{\bm{G}}_{3I-1}(\eta)\hat{\bm{G}}_{3I}(\zeta)$. Meanwhile, $F_I(\bm{\xi})$ =$(\bm{h}(\eta))^{\mathrm{T}}_{i,:}\bm{g}(\eta)^{\mathrm{T}}\bm{f}(\xi)^{\mathrm{T}}$. Thus, after passing the global random test, Eq.(\ref{C3ijk}) can be transformed to Eq.(\ref{C3ijka})
\begin{equation}\label{C3ijka}
\bm{\breve{C}}_3 \approx \bm{A}_{1} \times^1 \bm{A}_{2} \times^1 \bm{A}_{3}
\end{equation}
where $ \times^1 $ is the mode(3,1) contracted product and each $\bm{A}_i$ is a third-order tensor with each slice $\bm{A}_{i(:,j,:)}$ expressed in Eq.(\ref{slice}).
\begin{equation}\label{slice}
\bm{A}_{i(:,j,:)}= \int_0^1 (\bm{h}(\zeta))^{\mathrm{T}} \int_0^1 \bm{g}(\eta)^{\mathrm{T}} \int_0^1 \bm{f}(\xi)^{\mathrm{T}}\hat{\bm{G}}_{3i-2(j,:)}(\xi)\mathrm{d}\xi \hat{\bm{G}}_{3i-1}(\eta)\mathrm{d}\eta\hat{\bm{G}}_{3i}(\zeta)\mathrm{d}\zeta
\end{equation}
After obtaining $\bm{\breve{C}}_3$, it is possible to further reduce the dimension of $\bm{\gamma}(\theta)$ by applying an orthogonal transformation represented by HOSVD bases of $\breve{C}_3(i,j,k)$. Since cumulant tensors are super-symmetric in theory, the modes are identical in all dimensions and need to be computed in only a single dimension. To achieve this goal, first, we derive Eq.(\ref{C3ijkjk})
\begin{equation}\label{C3ijkjk}
\begin{split}
\bm{\breve{C}}_{3(1)}\bm{\breve{C}}_{3(1)}^{\mathrm{T}} &= \sum\limits_{j,k} (\bm{A}_{1(2)}\bm{A}_{2(:,j,:)}\bm{A}_{3(:,k,:)})(\bm{A}_{1(2)}\bm{A}_{2(:,j,:)}\bm{A}_{3(:,k,:)})^{\mathrm{T}}\\
&= \bm{A}_{1(2)}(\sum\limits_j \bm{A}_{2(:,j,:)} (\sum\limits_k \bm{A}_{3(:,k,:)}\bm{A}_{3(:,k,:)}^{\mathrm{T}})\bm{A}_{2(:,j,:)}^{\mathrm{T}})\bm{A}_{1(2)}^{\mathrm{T}}
\end{split}
\end{equation}
where subscripts $(i)$ is the $i$-mode unfolding of a tensor. It is clear in Eq.(\ref{C3ijkjk}) that the summation with respect to $j,k$ is decoupled. Then the pursued orthogonal transformation matrix $\bm{U}_{[3]}$=$ \bm{U}_{(:,1:n)} $ is obtained by spectral decomposition of matrix $\bm{\breve{C}}_{3(1)}\bm{\breve{C}}_{3(1)}^{\mathrm{T}}$. The number of columns $ n $ in $\bm{U}$ is selected as the minimum value satisfying Eq.(\ref{Un})
\begin{equation}\label{Un}
\min\left(\dfrac{\lVert(\bm{U}_{[3]}^{\mathrm{T}})\bm{\lambda}_2(\bm{U}_{[3]})\rVert_2}{\lVert\bm{\lambda}_2\rVert_2}, \dfrac{\lVert\bm{U}_{[3]}^{\mathrm{T}}\bm{\lambda}_3\bm{U}_{[3]}\rVert_2}{\lVert\bm{U}^{\mathrm{T}}\bm{\lambda}_3\bm{U}\rVert_2}  \right)>tol_3\ (0<tol_3<1)
\end{equation}
where $ \bm{\lambda}_2 $ and $ \bm{\lambda}_3 $ are diagonal matrices containing eigenvalues of $\tilde{C}_{-\xi\eta}(\zeta,\zeta')$ and mode-1 eigenvalues of $\bm{\breve{C}}_{3(1)}$, respectively. Next, based on the properties of TT decomposition \cite{Lee2018}, the cumulant tensor of the transformed latent factors $\bm{\gamma}_{[3]}(\theta) $ = $\bm{U}_{[3]}^{\mathrm{T}}\bm{\gamma}(\theta)$ is obtained by Eq.(\ref{C3new})
\begin{equation}\label{C3new}
\begin{split}
\bm{\acute{C}}_3& = \bm{\breve{C}}_3 \times_1 \bm{U}_{[3]}^{\mathrm{T}} \times_2 \bm{U}_{[3]}^{\mathrm{T}} \times_3 \bm{U}_{[3]}^{\mathrm{T}}\\
 &= \bm{\acute{A}}_{1} \times^1 \bm{\acute{A}}_{2}\times^1 \bm{\acute{A}}_{3}
\end{split}
\end{equation}
where $\bm{\acute{A}}_{j(:,:,i)}$ = $\bm{A}_{j(:,:,i)}\bm{U}$. The final expression of the random field is obtained as Eq.(\ref{RVfinal})
\begin{subequations}\label{RVfinal}
\begin{alignat}{1}
(x,y,z) &= \sum_{\bm{I}}R_{\bm{I}}^{\bm{p}}(\xi,\eta,\theta)\bm{B_I}\\
\begin{split}
\alpha(\xi,\eta, \zeta;\theta) &\approx (\bm{f}(\xi)_{1\times n_1} \bm{g}(\eta)_{n_1\times n_2}\bm{h}(\zeta)_{n_2\times n_3}\bm{U}_{[3]n_3\times n}) \bm{\gamma}_{[3]}(\theta)_{n\times 1}\\
&= \bm{F}(\xi,\eta,\zeta)\bm{\gamma}_{[3]}(\theta)
\end{split}\\
\mathrm{Cum}_1(\bm{\gamma}_{[3]}) &=\bm{0}\\
\mathrm{Cum}_2(\bm{\gamma}_{[3]}) &=\bm{U}_{[3]}^{\mathrm{T}}\bm{\lambda}_2 \bm{U}_{[3]}\\
\mathrm{Cum}_3(\bm{\gamma}_{[3]}) &=\bm{\acute{C}}_3
\end{alignat}
\end{subequations}
Accuracy of this expression for matching the prescribed second and third-order cumulant function is checked by computing the empirical global relative error in Eq.(\ref{epgf})
\begin{equation}\label{epgf}
\varepsilon_{\mathrm{gf},k} = \dfrac{\sqrt{\dfrac{1}{N}\sum\limits_{i=1}^{N}(\tilde{C}(\xi_i,\eta_i,\zeta_i; \xi'_i,\eta'_i,\zeta'_i) - \tilde{C}_{\mathrm{TTf},k}(\xi_i,\eta_i,\zeta_i; \xi'_i,\eta'_i,\zeta'_i))^2 }}{\sqrt{\dfrac{1}{N}\sum\limits_{i=1}^{N}\tilde{C}(\xi_i,\eta_i,\zeta_i; \xi'_i,\eta'_i,\zeta'_i)^2 }} 
\end{equation}
where $ \tilde{C}_{\mathrm{TTf},k} $ is the $ k $-th order cumulant function of the final TT expression and is computed by Eq.(\ref{CTTf})
\begin{subequations}\label{CTTf}
\begin{alignat}{1}
\tilde{C}_{\mathrm{TTf},2}(\xi_i,\eta_i,\zeta_i; \xi'_i,\eta'_i,\zeta'_i)&= \bm{F}(\xi_i,\eta_i,\zeta_i)(\bm{U}_{[3]}^{\mathrm{T}}\bm{\lambda}_2 \bm{U}_{[3]})\bm{F}(\xi_i,\eta_i,\zeta_i)^{\mathrm{T}}\\
\begin{split}
\tilde{C}_{\mathrm{TTf},3}(\xi_i,\eta_i,\zeta_i; \xi'_i,\eta'_i,\zeta'_i)&= \bm{\acute{C}}_3 \times_1 \bm{F}(\xi_i,\eta_i,\zeta_i) \times_2 \bm{F}(\xi_i,\eta_i,\zeta_i) \times_3 \bm{F}(\xi_i,\eta_i,\zeta_i)\\
&= (\bm{\acute{A}}_1 \times_2 \bm{F}(\xi_i,\eta_i,\zeta_i))\times^1 (\bm{\acute{A}}_2 \times_2 \bm{F}(\xi_i,\eta_i,\zeta_i))\\
&\quad\times^1  (\bm{\acute{A}}_3 \times_2 \bm{F}(\xi_i,\eta_i,\zeta_i))
\end{split}
\end{alignat}
\end{subequations}
If $ \min(\varepsilon_{\mathrm{gf}k}) $ is smaller than a prescribed threshold, then finish; else, decrease $ tol_3 $ and re-compute each $ \varepsilon_{\mathrm{gf}k} $. The above procedure is readily extended to other values of $m$. 
\subsection{Summary}
For the sake of clarity, the whole algorithm framework is summarized as follows:
\begin{enumerate}
	\item Do space transformation;
	\item Compute TT decomposition of $ \tilde{C} $ and perform global random test;
	\item Compute modes in each direction;
	\item Compute TT decomposition of $ \tilde{C}_3 $ and perform global random test;
	\item Compute $ \bm{\breve{C}}_3 $ and its HOSVD;
	\item Compress modes and latent factors, and perform the final global random test.
\end{enumerate}

When multiple higher-order cumulant functions are simultaneously given, we need to evaluate their relative importance and find an orthogonal transformation which can best explain a functional of these cumulant functions. This is a non-trivial task. The simple criterion given in \cite{Morton2009} lacks objectivity. Hence, this problem still needs further research.

\section{Comparing with related works}\label{Compa}

\subsection{Tensor train-PCA}

A recent work worth mentioning is the tensor train-PCA \cite{Wang2019}. Given $N$ samples of tensor data $\mathcal{X}_i\in\mathbb{R}^{I_1\times\cdots I_n}$, $i = 1,\dotsc, N$, the goal is to find $\mathcal{U}_j\in\mathbb{R}^{r_{j-1}\times I_j\times r_j}$, $j = 1,\dotsc, n$, such that the distance of the points to the space spanned by $\{\mathcal{U}_j\}_{j=1}^{n}$ is minimized. The algorithm in the original article requires concatenation of all the samples, hence is in essence a Full-to-TT algorithm \cite{Oseledets2010}, and is a data-driven realization of a part of our theoretical framework in Section \ref{theocore} as well. The modes are computed by applying SVD on the unfolded data matrices, which is computationally infeasible for large values of $N$ and $I_j$, $j = 1,\dotsc, n$. Different from the tensor train-PCA, our goal is to find modes of a \textit{continuously} indexed random field in each direction with \textit{prescribed} cumulant functions, as in Section \ref{theocore}. Estimating cumulant functions from raw data is out of the scope of this work. In addition, manipulations of large-scale arrays never occur in our work.

\subsection{Independent component analysis}

According to \cite{DeLathauwer2000,Nordhausen2018}, with the key assumption that the latent factors are mutually independent, independent component analysis (ICA) consists of the two sequential steps: (1) standardization (also called pre-whitening) of the original random field using the covariance function and (2) rotation of the standardized random field to the latent factors named independent components by using higher-order (usually the fourth) cumulant functions. Thus, the proposed theoretical framework has a similar form as that of ICA. However, rather than pursuing the independent latent factors, the proposed framework just aims at reducing the dimensionality of the latent factors and represent their higher-order non-Gaussian dependence.

\subsection{Isogeometric analysis-based K-L expansion}

IGA-based K-L expansion were proposed in recent papers to represent random fields on complex domains. However, the underlying architecture of these methods is still trivial PCA in \textit{physical} domain. Modes are projected onto a subspace spanned by IGA basis functions, and the coordinates are computed with Galerkin or collocation method. Hence, these methods still suffer from the drawbacks of trivial PCA. Different from these methods, the underlying architecture of this work is the \textit{tensor train} decomposition in \textit{parametric} domain. Highlights of the proposed framework over the IGA-based methods is summarized as follows:
\begin{enumerate}
	\item Since we only need a space transformation and the geometry is exactly preserved at the lowest level of refinement, enrichment of NURBS spaces is never needed.
	\item The original \textit{unstructured} random field is transformed into a \textit{structured} one, making it feasible for using the powerful tensor train decomposition to compress the random field. The tensor train representation has the merits of both the trivial PCA (low stochastic dimensionality) and other tensor-based representations (moderate physical dimensionality).
	\item It is totally unnecessary to predefine any basis functions or discrete point set to represent the modes, which avoids subjectivity. All modes are automatically represented in \emph{chebfun} format in the sense of machine precision.
	\item By adaptive variable separation of the covariance function, all integrals are transformed into uni-dimensional ones, which greatly reduces the computational scale of each mode. In addition, the tedious processes of forming and assembling the stiffness matrices are never needed. All the eigenvalues are obtained by solving only a few (1-3) small-scale ($ O(10^{0-2}) $) matrix SVD problems.
	\item Since the random field is represented in the parametric domain rather than the physical counterpart, Jacobians are never needed in the integrals. 
	\item By adaptive variable separation of the random field and the higher-order cumulant functions, cumulants of latent factors can be conveniently computed to capture the higher-order non-Gaussian correlations of the random field, which has never been discussed in existing literature of K-L expansion.
\end{enumerate}
\textbf{Remark}: Since the random field is represented in different spaces, results of the proposed algorithm framework are not directly comparable with that of existing IGA-based methods.

\section{Case studies}\label{Cases}

In this section, we examine the performance of the proposed algorithms by using three examples with increasing parametric dimensionalities. The underlying isogeometric transformation is achieved by using the toolbox in \cite{Vazquez2016}. All the numerical experiments were performed using MATLAB (Version 2019a) on a notepad (core i5 CPU and 12GB RAM).

\subsection{Example 1}

In the first example, to clearly demonstrate performance of the proposed method, we consider a simple random field with $D=[0,1]\subseteq\mathbb{R}$ and covariance function defined in Eq.(\ref{ex1Cum2})
\begin{equation}\label{ex1Cum2}
C(x,x')=\sum\limits_{k=1}^{80} \lambda_k f_k(x)f_k(x') 
\end{equation} 
where $\lambda_k = 4/(\pi^2(2k-1)^2)$ and $f_k(x) = \sqrt{2}\sin(\xi/\sqrt{\lambda_k})$. Following the steps in Section \ref{Compu}, first, the isogeometric transformation is constructed with the following ingredients: (1) knot vector $\Xi$ = (0,0,1,1); (2) control points $x_1$ =0, $x_2$ =1; (3) weights $\bm{W}$ = (1,1). It is easy to derive that the transformation is $x = \xi$. Then, eigenpairs of the compound integral operator with kernel $C_2(x(\xi),x(\xi'))$ should be computed. Since the parametric space coincide with the physical space, the eigenpairs are exactly $ \lambda_k $ and $f_k(\xi)  $. Thus, TT-rank of $ \tilde{C} $ is exactly 80. Comparisons are made from different aspects as below.

\subsubsection{Visualization of the reconstructed compound covariance function and  spatial distribution of the errors}

\begin{figure}
	\centering
	\subfigure[reconstructed compound kernel]{\includegraphics[width =0.4\textwidth]{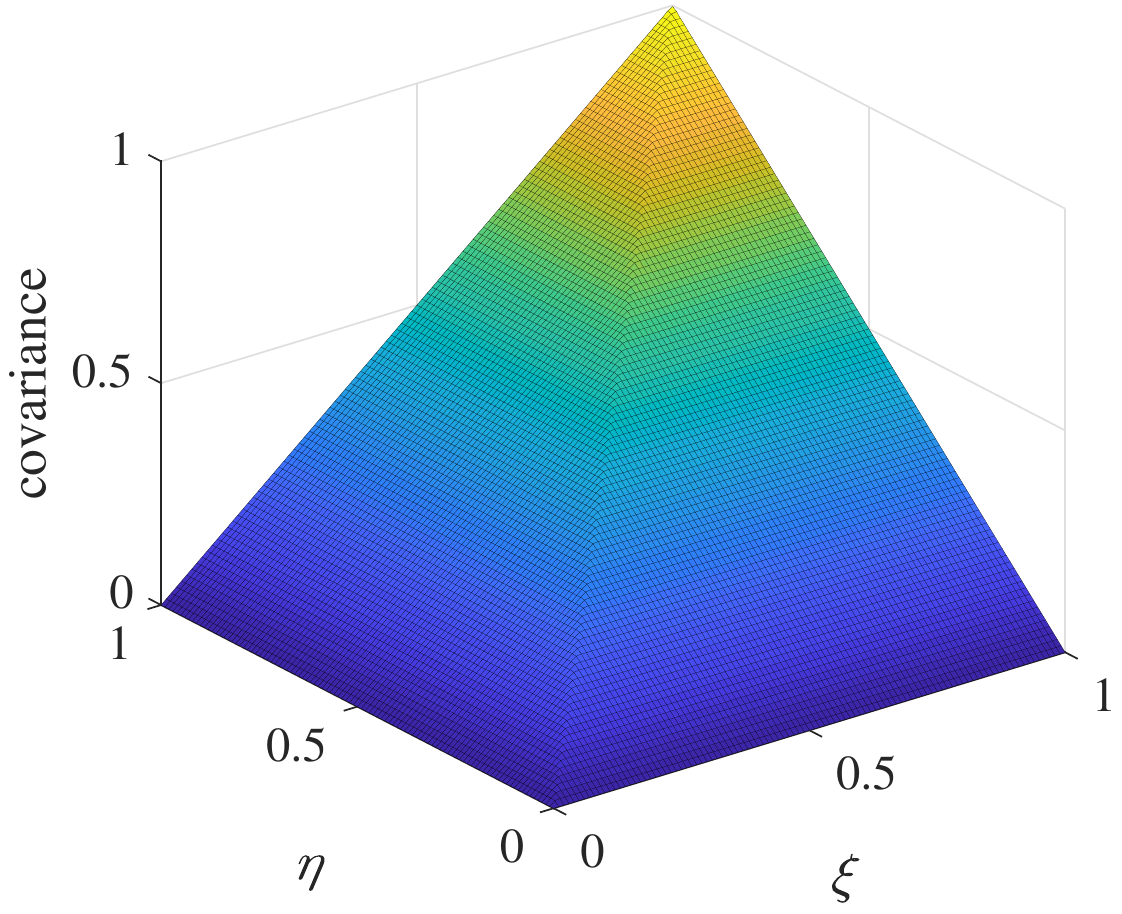}\label{ex1-1-1}}	
	\subfigure[spatial distribution of the errors]{\includegraphics[width = 0.4\textwidth]{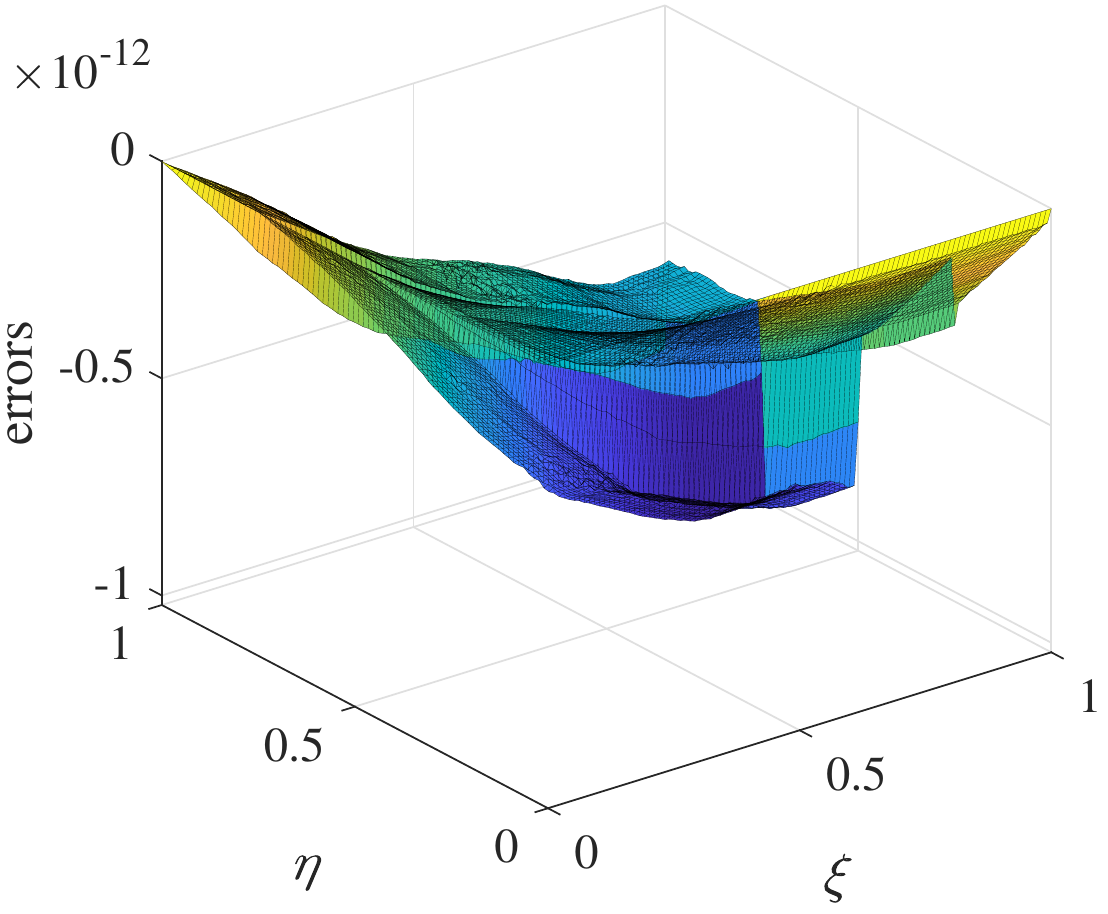}\label{ex1-1-2}}	
	\subfigure[the interpolation set]{\includegraphics[width = 0.4\textwidth]{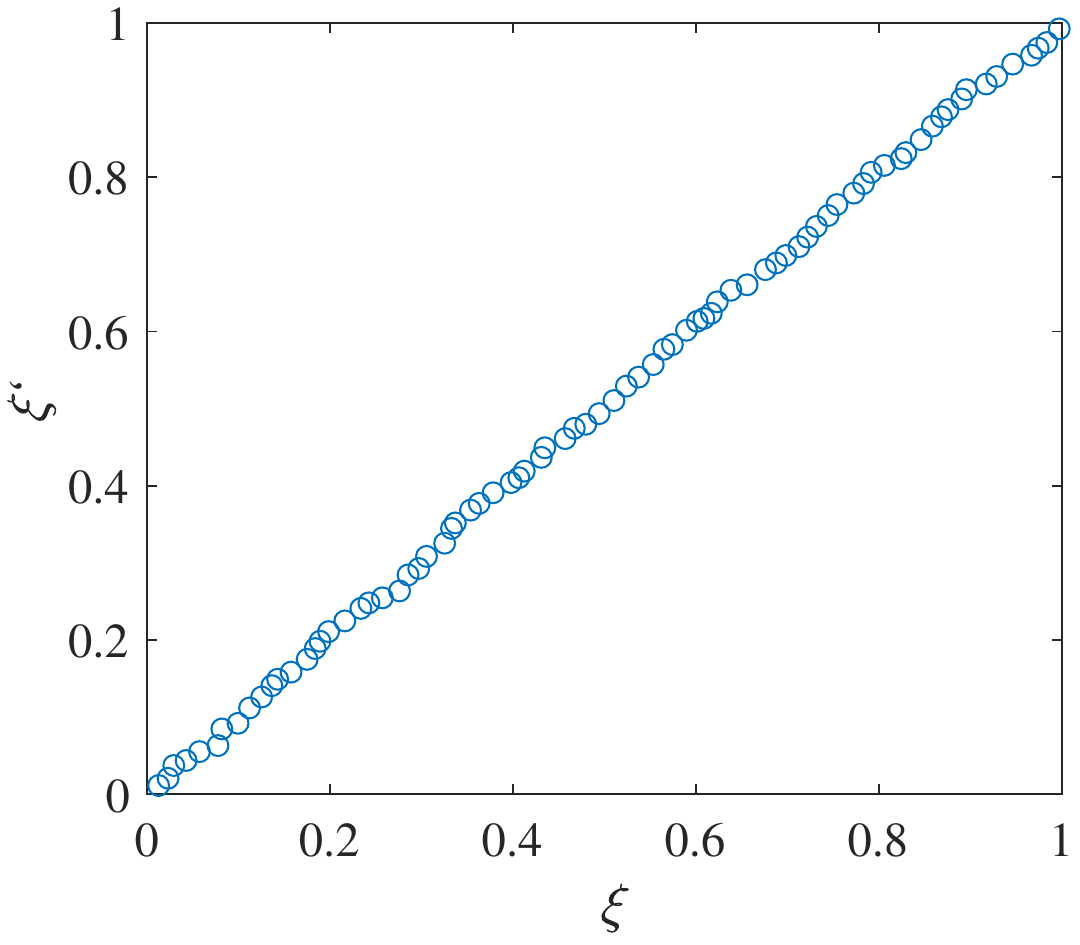}\label{ex1-1-3}}
	\caption{Visualization of the results of a numerical experiment in Example 1}\label{ex1-1}	
\end{figure}

First, an experiment is made by letting $N$=400 and $tol=1.0000\times 10^{-6}$. Results are illustrated in Fig.\ref{ex1-1}.

Fig.\ref{ex1-1-1} shows the reconstructed covariance function. Spatial distribution of errors in Fig.\ref{ex1-1-2} indicates that this reconstruction is almost perfect since errors are consistently of $ O(10^{-12}) $. In addition, the interpolation set is automatically selected as the points shown in Fig.\ref{ex1-1-3}, which indicates the abilities of the proposed algorithm for detecting the regularity of the target function and adapting to its features.

\subsubsection{Investigations on global relative errors and ranks} 

\begin{figure}
	\centering
	\subfigure[global relative error]{\includegraphics[width = 0.4\textwidth]{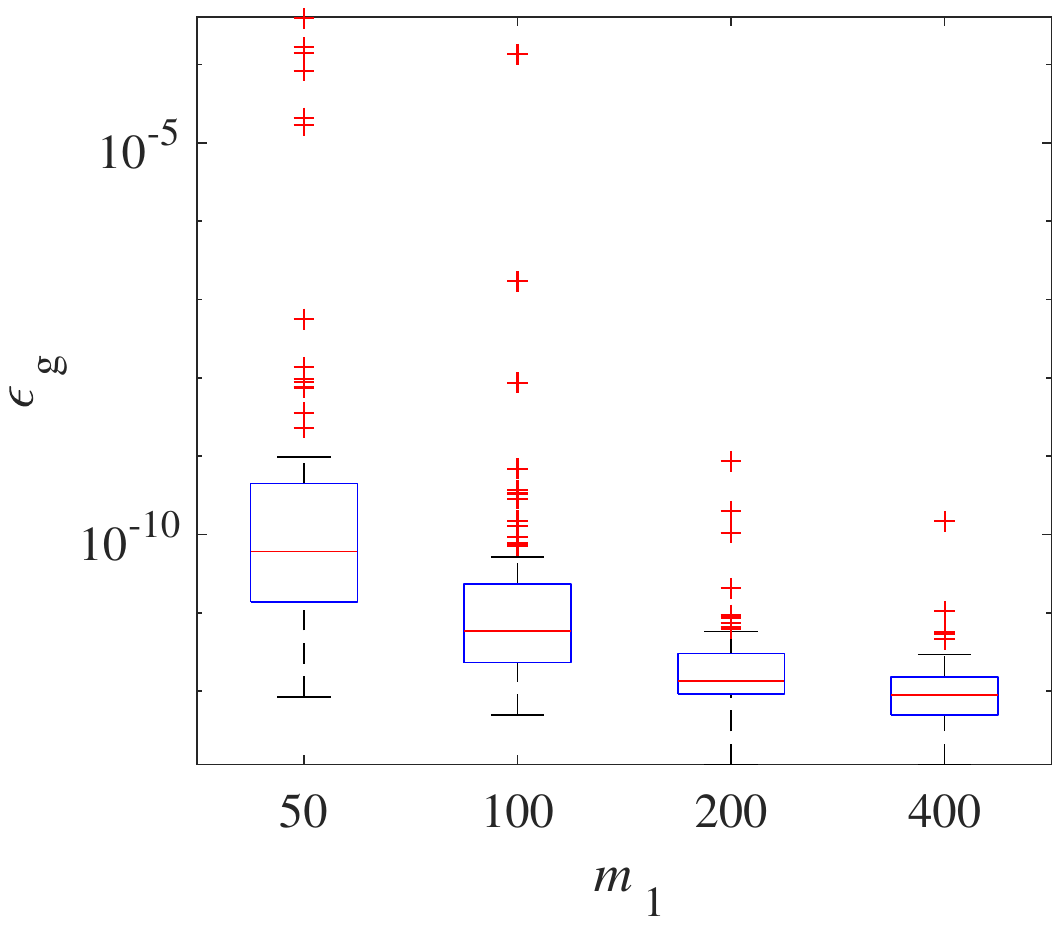}\label{ex1-2-1}}		
	\subfigure[rank]{\includegraphics[width = 0.4\textwidth]{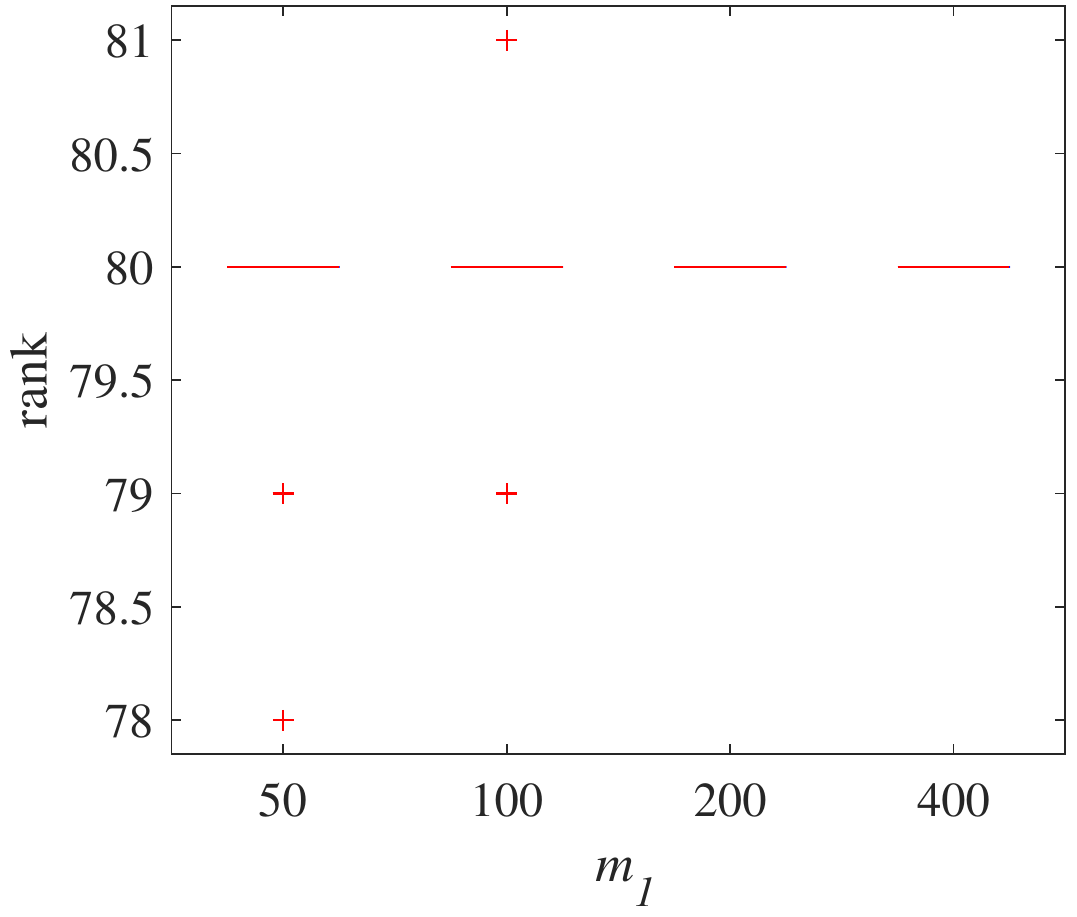}\label{ex1-2-2}}
	\caption{Variability of the results with different $m_1$ values}\label{ex1-2}	
\end{figure}

 To investigate the effect of $m_1$ in line \ref{m1} of Algorithm \ref{ctt2} on variability of the global relative errors and the rank, we did four groups of experiment with $m_1$ = 50, 100, 200 and 400, respectively. Each group consisted of 100 independent tests. Each test was performed by letting $N$ = 1000 in Eq.(\ref{epg}) and $tol$ = $10^{-6}$ in Algorithm \ref{ctt2}. Results are illustrated in Figs.\ref{ex1-2}.
Fig.\ref{ex1-2-1} shows the number of outliers decreases from 14 to 5 as $ m_1 $ increases from 50 to 400. Despite the outliers, both the variability (represented by difference between the upper and lower adjacent value) and the median of  $\varepsilon_{\mathrm{g}}$ decreases with $ m_1 $ as well. Moreover, the decreasing rates decreases with $m_1$. These three observations indicate that both the average accuracy and the robustness of the proposed algorithm increases with $ m_1 $, while the increasing rates tend to be slower. When $ m_1=400 $, very accurate (median=$ O(10^{-12}) $) and robust (upper adjacent value=$ O(10^{-12}) $ and rank $ \equiv $ 80) can be obtained. Meanwhile, Fig.\ref{ex1-2-2} shows that the proposed algorithm almost always reveals the exact rank (except only a few outliers) with all values of $ m_1 $.

\begin{figure}
	\centering
	\subfigure[global relative error]{\includegraphics[width = 0.4\textwidth]{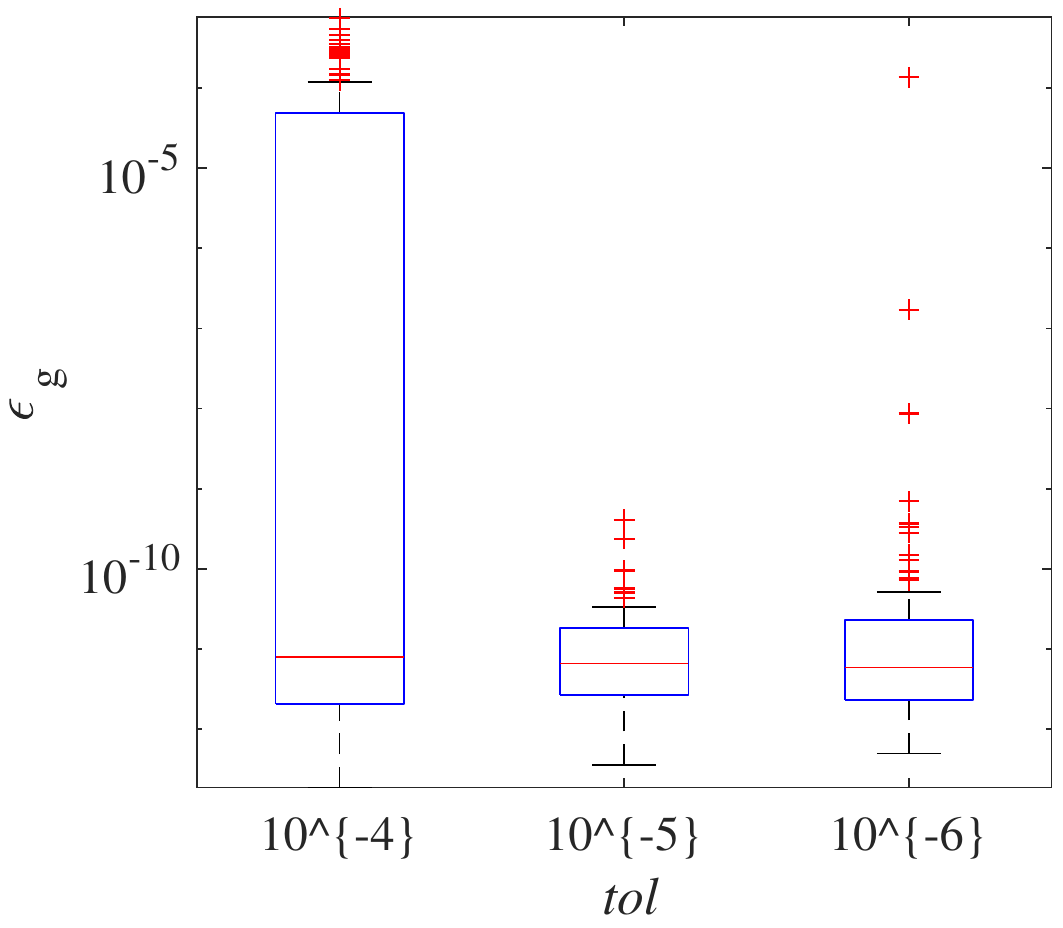}\label{ex1-3-1}}	
	\subfigure[rank]{\includegraphics[width = 0.4\textwidth]{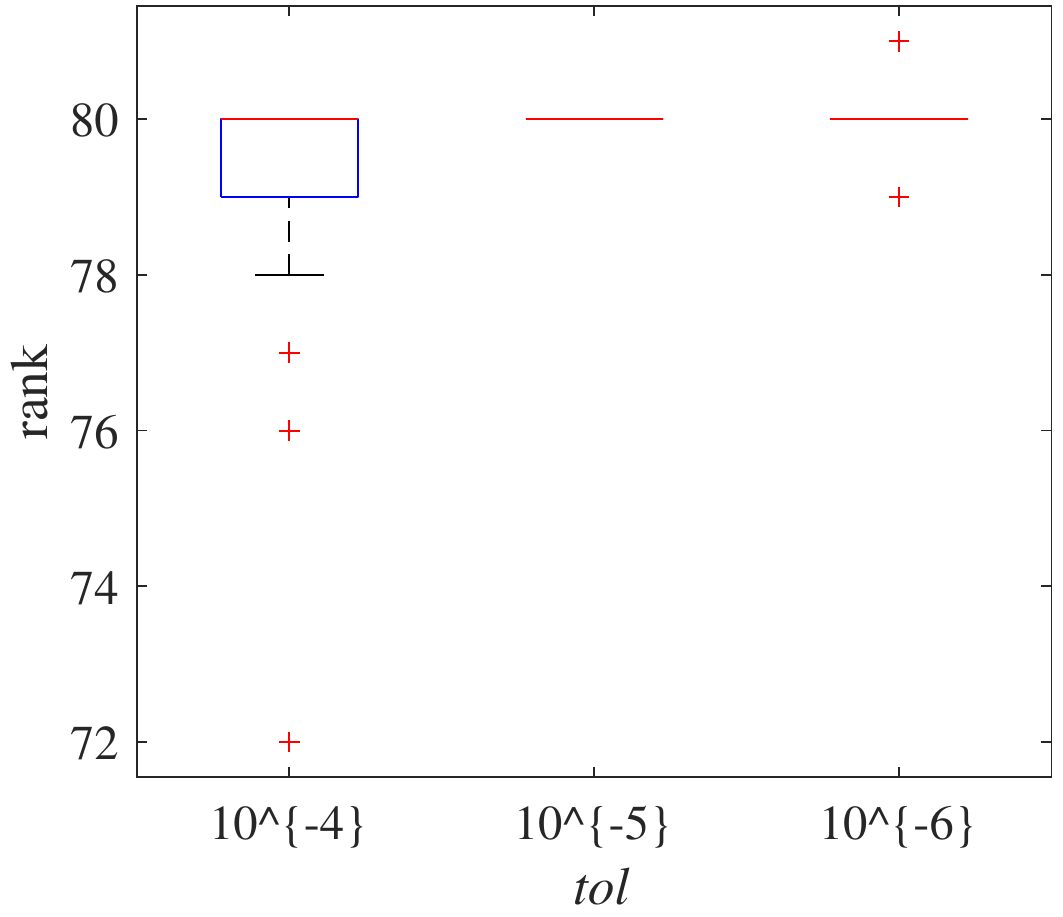}\label{ex1-3-2}}
	\caption{Variability of the results with different $tol$ values}\label{ex1-3}	
\end{figure}

Next, two additional groups of experiment were made to investigate the effects of $tol$ on the variability of the global relative errors and the rank by setting $tol$ = $10^{-4}$ and $10^{-5}$ , respectively. $m_1$ is kept as 100 in both groups. Other settings are the same as that of the previous four groups. Results are illustrated in Fig.\ref{ex1-3}.
Fig.\ref{ex1-3-1} shows that variability of $\varepsilon_{\mathrm{g}}$ significantly decreases as $ tol $ decreases from $ 10^{-4} $ to $ 10^{-5} $, then almost remains constant as $ tol $ decreased to $ 10^{-6} $. The median almost remains constant in all cases. Similar features are displayed in the figure of rank. Thus, to obtain sufficient robust results, $ tol $ is suggested to be of $ O(10^{-5}) $ or smaller in this example.

\begin{figure}
	\centering
	\subfigure[the 1st mode]{\includegraphics[width = 0.4\textwidth]{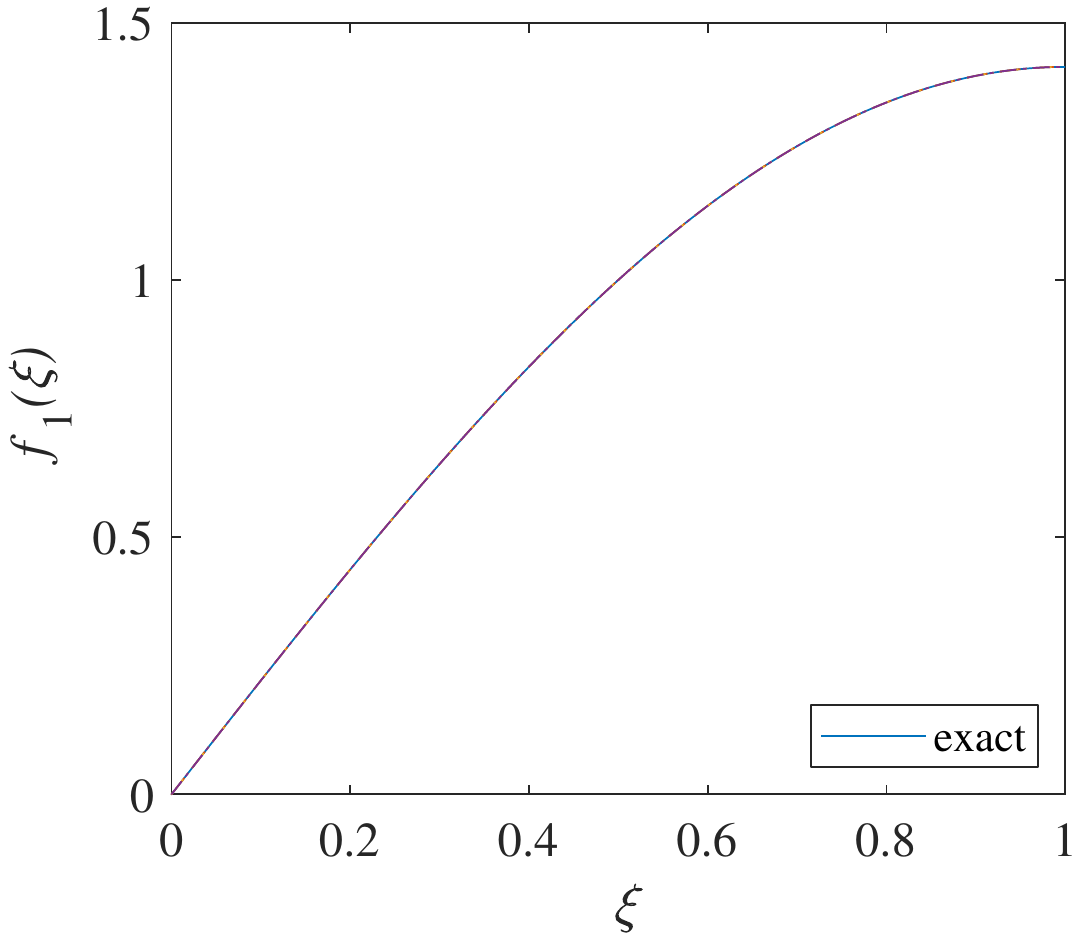}}	
	\subfigure[the 10th mode]{\includegraphics[width = 0.4\textwidth]{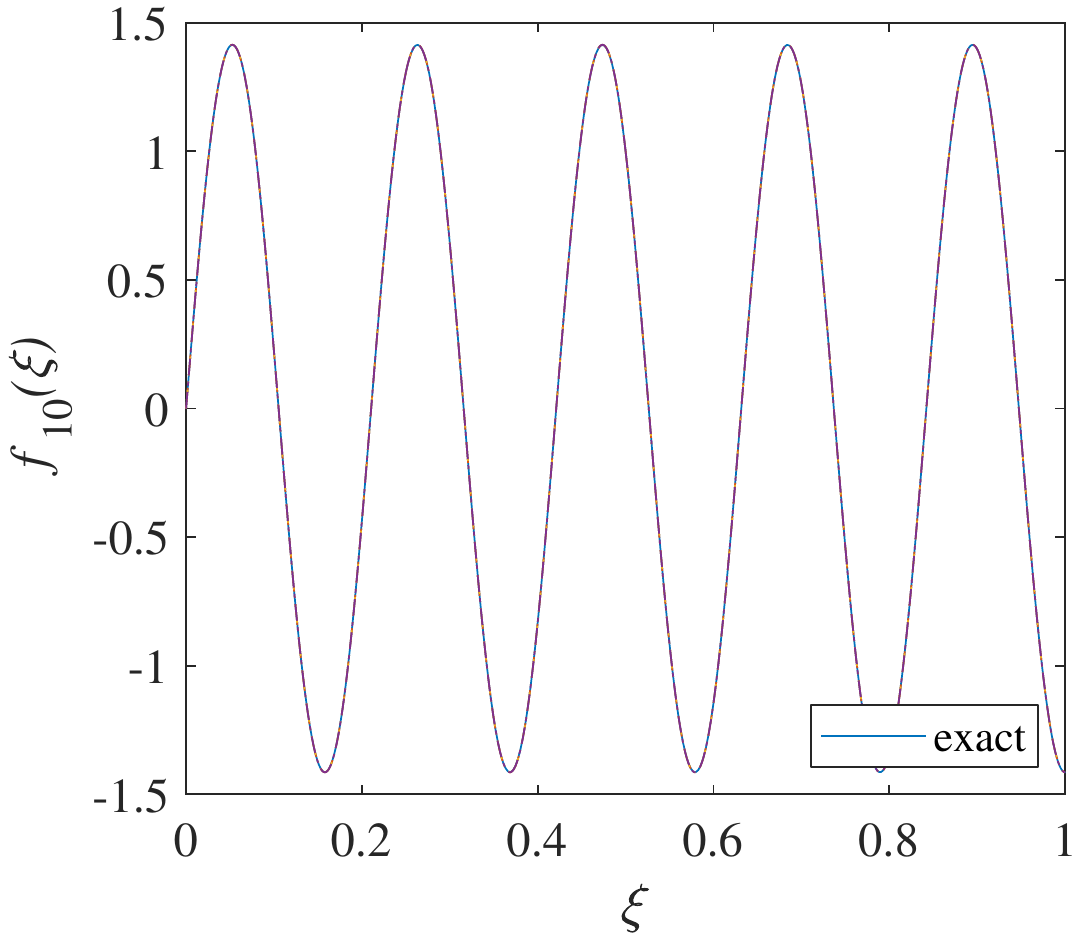}}	
	\subfigure[the 80th mode]{\includegraphics[width = 0.4\textwidth]{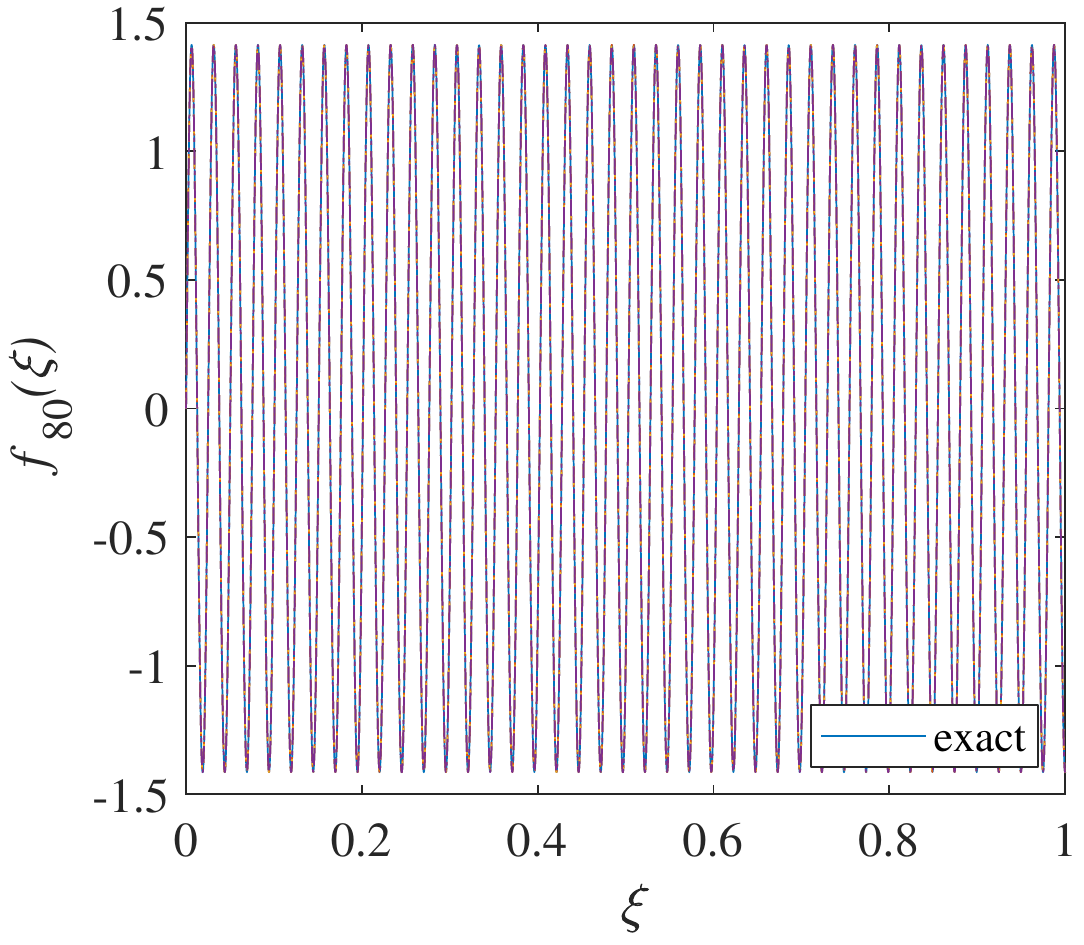}}	
	\subfigure[the 80 largest eigenvalues]{\includegraphics[width = 0.4\textwidth]{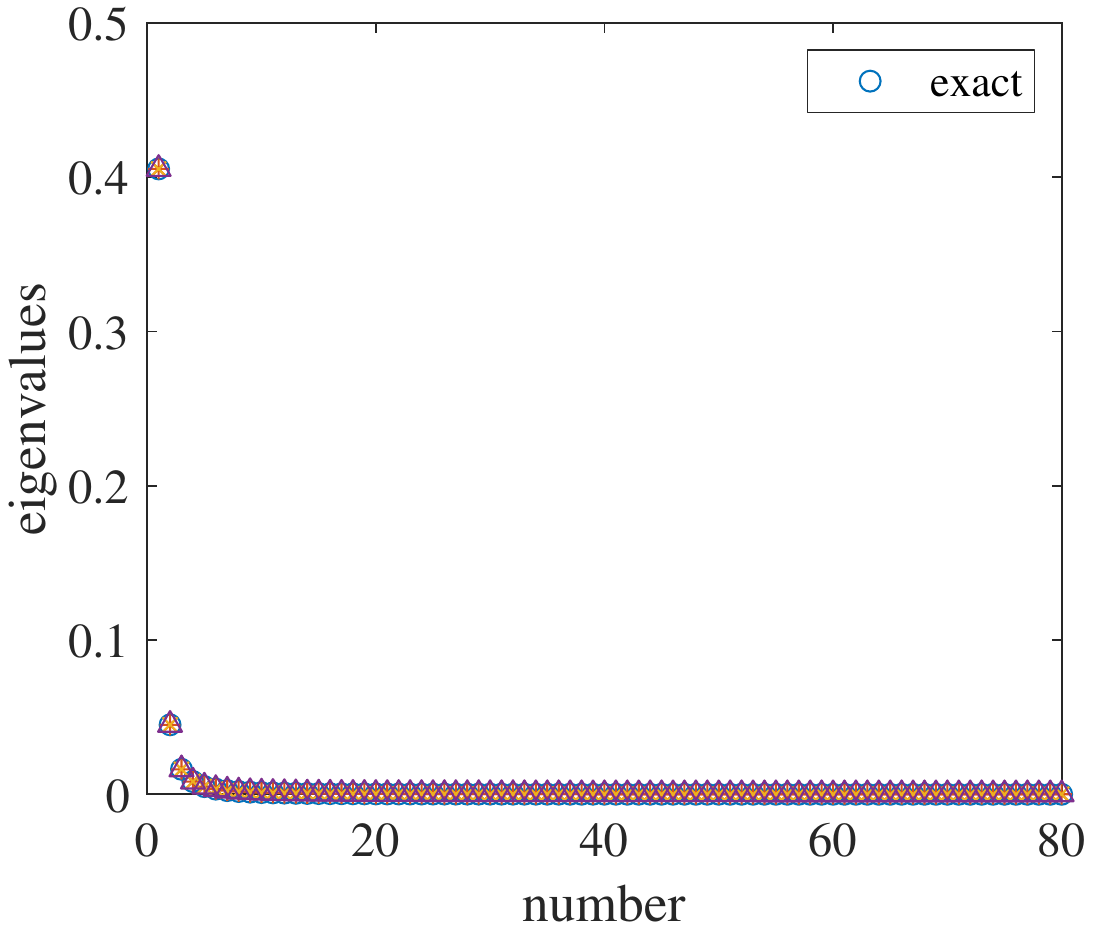}}
	\caption{Comparison of modes and eigenvalues in Example 1}\label{ex1-4}	
\end{figure}

Next, three independent tests were proceeded by setting $m_1$ = 400 and $tol$ = $10^{-6}$ to illustrate the accuracy and robustness of Algorithm \ref{ctt2} for computing the modes. Results are shown in Fig.\ref{ex1-4}.
Fig.\ref{ex1-4} shows that for all orders of mode (from 1 to 80), the proposed algorithm consistently provides almost perfect reconstructions (the differences can not be recognized by naked eyes) in all the three tests, and so do the eigenvalues. These observations verify the high levels of accuracy and robustness of the proposed method from another point of view.

\subsubsection{Dimension reduction of the latent factors} 

Before proceeding dimension reduction, we added an experiment with $K=2$, $m_1=400$ and $tol=10^{-6}$, and computed $\varepsilon_{\mathrm{g}2}$ = $3.8889\times 10^{-13}\ (N=1000)$. Then, rank-80 eigenpairs were computed with $\varepsilon_{\mathrm{gf}2}$ = $1.2715\times 10^{-2}\ (N=1000)$. Next, we assume the third-order cumulant function $C_3(x,y,z)$ as Eq.(\ref{ex1Cum3})
\begin{equation}\label{ex1Cum3}
C_3(x,y,z)=\sum\limits_{k=1}^{80} \lambda_k f_k(x)f_k(y)f_k(z) 
\end{equation} 
 and computed the TT decomposition $\tilde{C}_{3,\mathrm{TT}}$ with $\varepsilon_{\mathrm{g}3}$ = $1.1213\times 10^{-12}\ (N=1000)$ and ranks $\bm{r}$=$ (80,80) $. Finally, by setting $tol_3=9.9990\times 10^{-1}$, the effective dimensionality of the latent factors is greatly reduced from 80 to 5 with $\varepsilon_{\mathrm{gf}3}$ = $1.2555\times 10^{-2}\ (N=1000)$. The time cost (seconds) of each quantity is listed in Table \ref{eg1tb1}.
\begin{table}[htbp]
	\centering
	\caption{Time cost (seconds) of the quantities in Example 1}\label{eg1tb1}
	\begin{tabular}{*{6}{c}}
		\toprule
		&$ \tilde{C}_{k,\mathrm{TT}} $ & $\varepsilon_{\mathrm{g}k}$ & modes & $\varepsilon_{\mathrm{gf}k}$&total\\
		\midrule
		$ k =2$  & 107   & 1        & 440   & 1 & 549 \\ 
		$ k =3$   & 4120 & 2 & 1467 & 1 & 5590 \\
		\bottomrule
	\end{tabular}
\end{table}
The results above indicate that the goal of matching the second and third-order cumulant functions within moderate time has been achieved.

\subsection{Example 2}

\begin{figure}
	\centering
	\includegraphics[width = 0.4\textwidth]{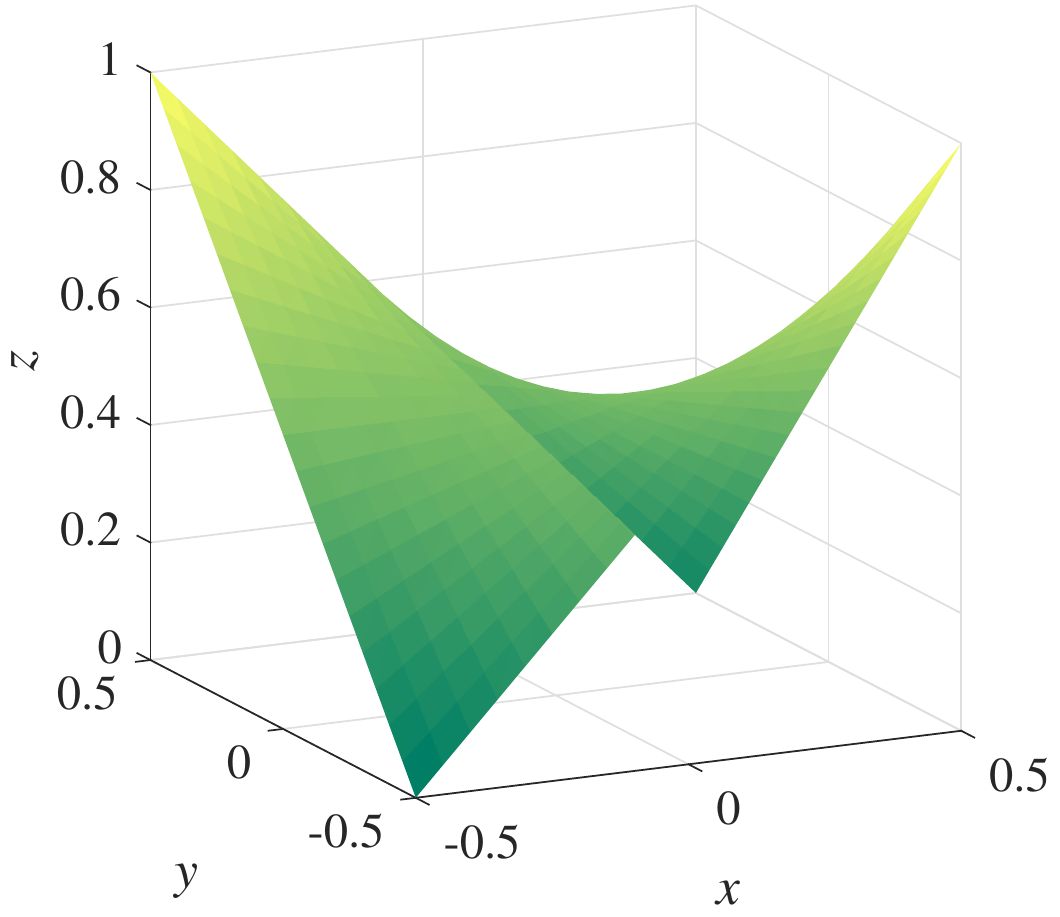}
	\caption{The physical domain in Example 2}\label{ex2-1}	
\end{figure}

In this example, let the physical domain $D$ be a bilinear surface in Fig.\ref{ex2-1} 
where the knot vectors are $\Xi_1$ = (0,0,1,1) and $\Xi_2$ = (0,0,1,1), and the control points and weights are listed in Table \ref{tbeg2}.
\begin{table}[htbp]
	\centering
	\caption{Control points and weights defining the geometry in Example 2}\label{tbeg2}
	\begin{tabular}{*{6}{c}}
		\toprule
		\multicolumn{2}{c}{Subscripts $\bm{i}$} & \multicolumn{2}{c}{Control points $B_{\bm{i}}$} & \multicolumn{2}{c}{weights $w_{\bm{i}}$} \\
		\midrule
		(1,1)  & (1,2) & (-0.5,-0.5,0) & (0.5,-0.5,1) & 1& 1 \\
		(2,1)  & (2,2) & (-0.5,0.5,1) & (0.5,0.5,0) & 1 & 1 \\
		\bottomrule
	\end{tabular}
\end{table}
The covariance function in physical domain is defined as Eq.(\ref{Coveg2})
\begin{equation}\label{Coveg2}
C(\bm{x},\bm{y})=\sigma^2 \exp\left(-\dfrac{\lVert \bm{x} - \bm{y} \rVert_2^2}{(bL)^2} \right)\quad \bm{x},\bm{y}\in D\subseteq\mathbb{R}^3
\end{equation}
where $\sigma^2=1$, $b=1$ and $L=1$.

\begin{figure}
	\centering
	\subfigure[the 1st mode]{\includegraphics[width = 0.45\textwidth]{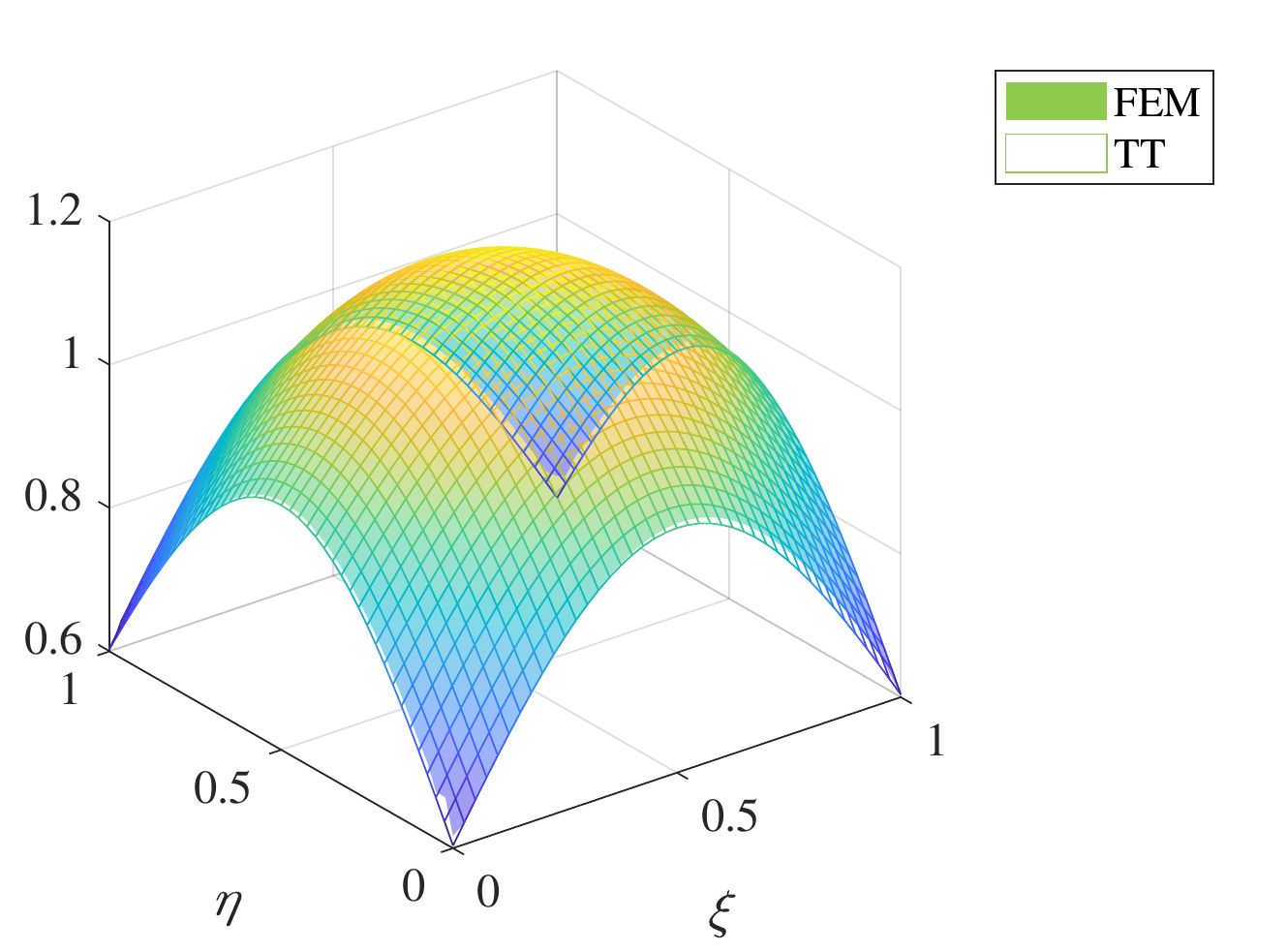}}
	\subfigure[the 2nd mode]{\includegraphics[width = 0.45\textwidth]{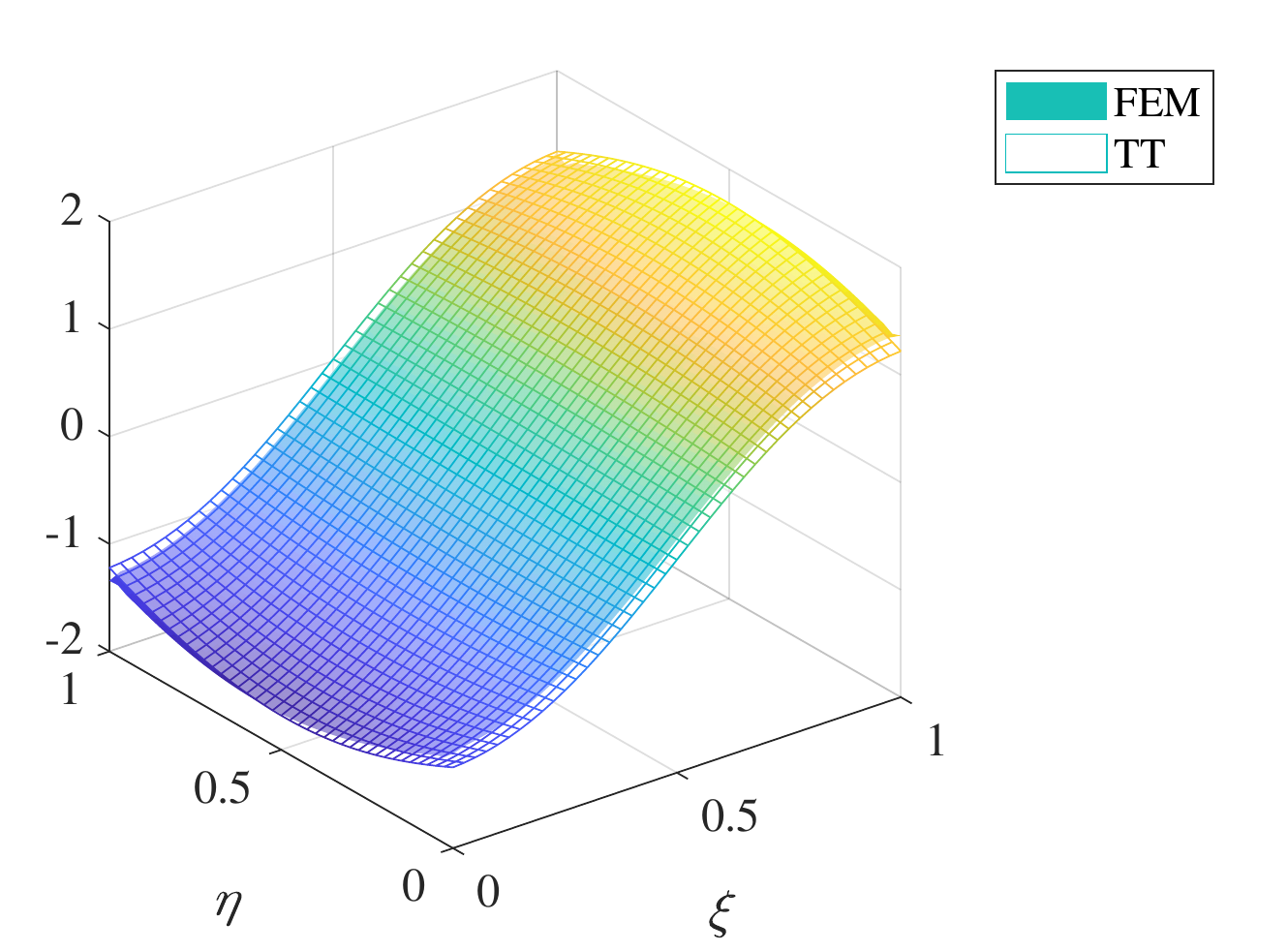}}
	\subfigure[the 4th mode]{\includegraphics[width = 0.45\textwidth]{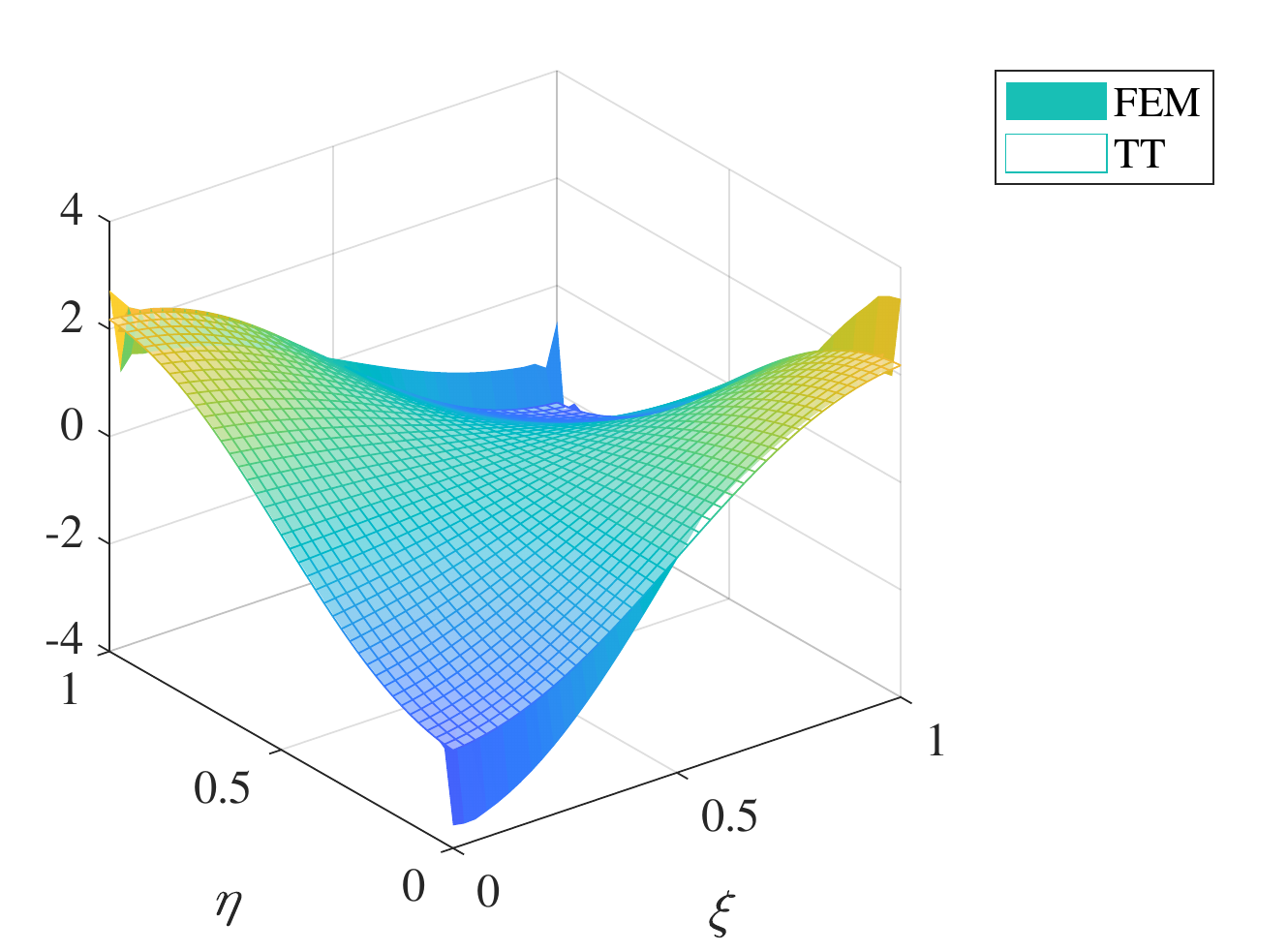}}
	\subfigure[the 37 largest eigenvalues]{\includegraphics[width = 0.45\textwidth]{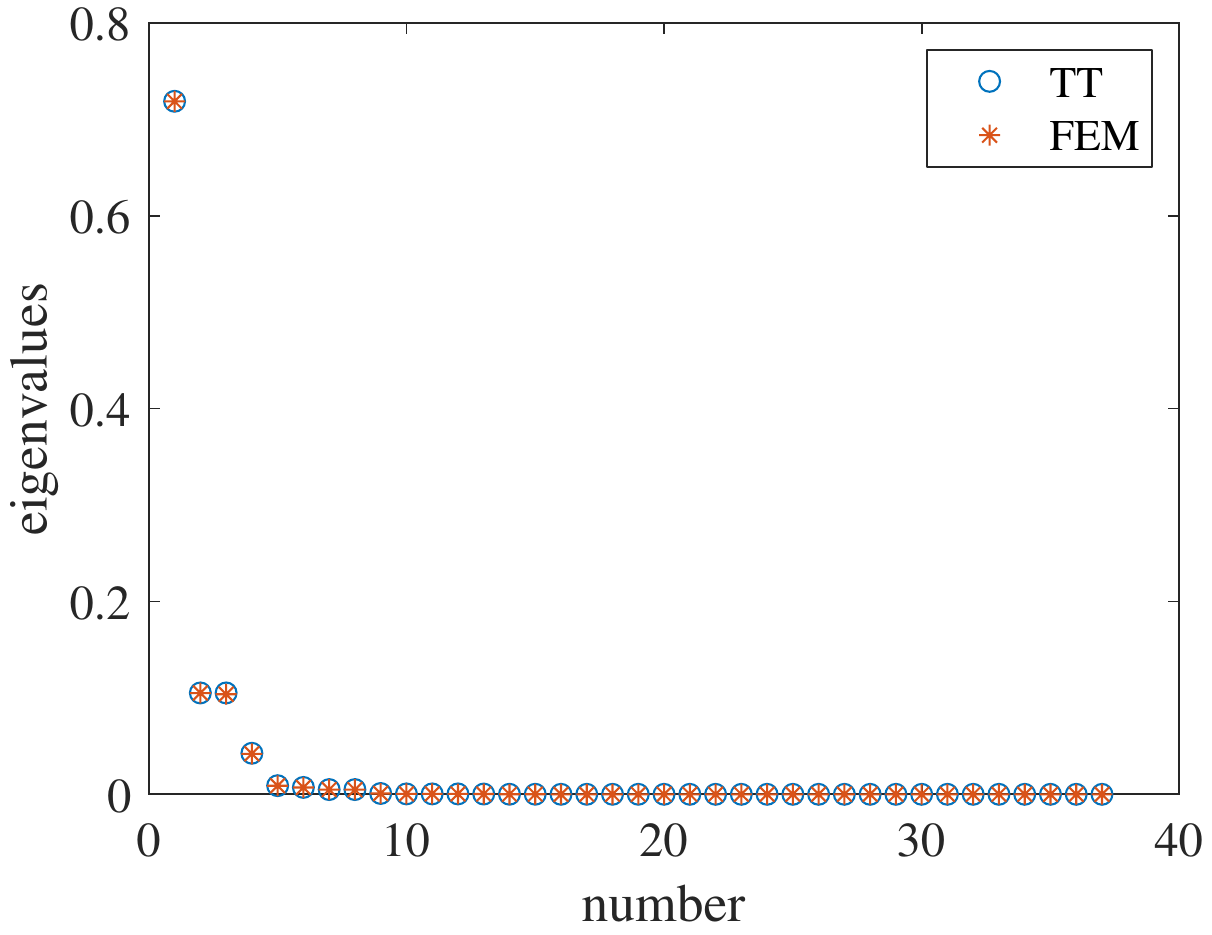}}
	\caption{Comparison of modes and eigenvalues in Example 2}\label{ex2-2}	
\end{figure}

First, we decomposed the covariance function in the parametric space by using Algorithm \ref{ctt} with $tol=1.0000\times 10^{-6}$ and each $m_k=800$ in Algorithm \ref{ISE}, and computed the global relative errors $\varepsilon_{\mathrm{g}2}=2.1065\times 10^{-7}(N=1000)$ and ranks $\bm{r}_2=(8,37,8)$. Then, we computed the eigenpairs in each direction with $\varepsilon_{\mathrm{gf}2}$ = $1.0847\times 10^{-2}\ (N=1000)$. To make comparisons, we also proceeded trivial PCA in the parametric space by using the FEM with rectangular bilinear elements and 40$ \times $40 mesh. Main results are illustrated in Fig.\ref{ex2-2}.

From Fig.\ref{ex2-2} we can see that the eigenvalues computed with the proposed method are highly consistent with the referenced values (the differences cannot be recognized by naked eyes). The dominant modes (whose corresponding eigenvalues are significantly larger than zero, the first four modes in this example) computed with the proposed method also show good agreement with the reference ones. Notice that the 'boundary effect' arise in the FEM solutions which themselves are approximate ones, accuracy of the TT modes cannot be evaluated merely by the discrepancies from the FEM solutions. Moreover, the TT solution is more smooth near the boundary since the 'boundary effect' hardly occurs. Thus, the TT solutions seem to be more accurate than the ones of FEM. On the whole, the above results are consistent with Theorem \ref{theo1}. 

Next, by assuming the third-order cumulant function as Eq.(\ref{C3eg2})
\begin{equation}\label{C3eg2}
C_3(\bm{x},\bm{y},\bm{z})=\sigma^3 \exp\left(-\dfrac{\lVert \bm{x} - \bm{y} \rVert_2^2 + \lVert \bm{x} - \bm{z} \rVert_2^2+\lVert \bm{y} - \bm{z} \rVert_2^2}{(bL)^2} \right)\quad \bm{x},\bm{y}\in D\subseteq\mathbb{R}^3
\end{equation}
where $\sigma^3=1$, $b=1$ and $L=1$, we computed the TT decomposition of the counterpart in the parametric domain (denoted as $ \tilde{C}_3(\xi,\eta;\xi',\eta';\xi'',\eta'') $) by setting $tol=1.0000\times 10^{-5}$ and each $m_k=800$ in Algorithm \ref{ISE}, and computed the global relative errors $\varepsilon_{\mathrm{g}3}=6.9247\times 10^{-6}(N=10000)$ and ranks $\bm{r}_3=(8,41,129,38,8)$. After HOSVD of $\bm{\breve{C}}_3$ ($tol_3=9.9990\times 10^{-1}$), the number of latent factors are reduced from 37 to 11 $ \varepsilon_{\mathrm{gf}3}=3.9148\times 10^{-3}, N=1000 $. The time cost (seconds) of each quantity is listed in Table \ref{eg2tb2}.
\begin{table}[htbp]
	\centering
	\caption{Time cost (seconds) of the quantities in Example 2}\label{eg2tb2}
	\begin{tabular}{*{6}{c}}
		\toprule
		 &$ \tilde{C}_{k,\mathrm{TT}} $ & $\varepsilon_{\mathrm{g}k}$ & modes & $\varepsilon_{\mathrm{gf}k}$&total\\
		\midrule
		$ k =2$   & 31   & 1        & 293   & 1 & 326 \\ 
		$ k =3$   & 575 & 3 & 3962 & 1 & 4540 \\
		\bottomrule
	\end{tabular}
\end{table}
 Computing the FEM solution above costs 5523s, which indicates that efficiency of the proposed algorithm is an order of magnitude higher than that of FEM. Thus, the goal of matching the second and third-order cumulant functions within moderate time has been achieved.

\subsection{Example 3}

\begin{figure}
	\centering
	\includegraphics[width = 0.4\textwidth]{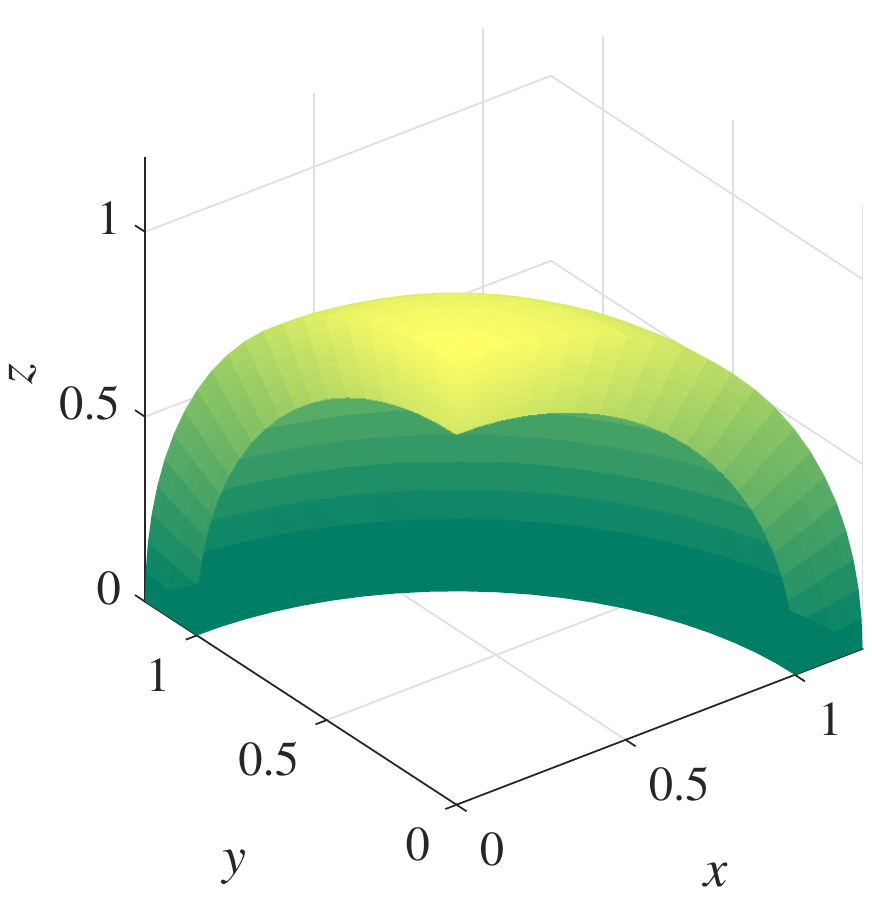}
	\caption{The physical domain in Example 3}\label{ex3-1}	
\end{figure}

In this example, let the physical domain $D$ be a hemispherical shell in Fig.\ref{ex3-1} 
where the inner radius $R_{\mathrm{i}}=1$ and the outer radius $R_{\mathrm{o}}=1.2$. The knot vectors are  $\Xi_1$ = (0,0,1,1), $\Xi_2$ = (0,0,0,1,1,1) and $\Xi_3$ = (0,0,0,1,1,1), respectively. The control points and weights are listed in Table \ref{tbeg3}.
\begin{table}[htbp]
	\centering
	\caption{Control points and weights defining the geometry in Example 3}\label{tbeg3}
	\begin{tabular}{*{6}{c}}
		\toprule
		\multicolumn{2}{c}{Subscripts $\bm{i}$} & \multicolumn{2}{c}{Control points $B_{\bm{i}}$} & \multicolumn{2}{c}{weights $w_{\bm{i}}$} \\
		\midrule
		(1,1,1) & (2,1,1) & $ (R_{\mathrm{i}},0,0) $  & $ (R_{\mathrm{o}},0,0) $ &  1& 1 \\
		(1,2,1) & (2,2,1) & $ (R_{\mathrm{i}},R_{\mathrm{i}},0) $  & $ (R_{\mathrm{o}},R_{\mathrm{o}},0) $  & $ \frac{1}{\sqrt{2}} $& $ \frac{1}{\sqrt{2}} $ \\
		(1,3,1) & (2,3,1) & $ (0,R_{\mathrm{i}},0) $  & $ (0,R_{\mathrm{o}},0) $ & 1& 1 \\
		
		(1,1,2) & (2,1,2) & $ (R_{\mathrm{i}},0,R_{\mathrm{i}}) $  & $ (R_{\mathrm{o}},0,R_{\mathrm{o}}) $  & $ \frac{1}{\sqrt{2}} $& $ \frac{1}{\sqrt{2}} $ \\
		(1,2,2) & (2,2,2) & $ (R_{\mathrm{i}},R_{\mathrm{i}},R_{\mathrm{i}}) $  & $ (R_{\mathrm{o}},R_{\mathrm{o}},R_{\mathrm{o}}) $  & $ \frac{1}{2} $& $ \frac{1}{2} $ \\
		(1,3,2) & (2,3,2) & $ (0,R_{\mathrm{i}},R_{\mathrm{i}}) $  & $ (0,R_{\mathrm{o}},R_{\mathrm{o}}) $  & $ \frac{1}{\sqrt{2}} $& $ \frac{1}{\sqrt{2}} $ \\
		
		(1,1,3) & (2,1,3) & $ (0,0,R_{\mathrm{i}}) $  & $ (0,0,R_{\mathrm{o}}) $ & 1& 1 \\
		(1,2,3) & (2,2,3) & $ (0,0,R_{\mathrm{i}}) $  & $ (0,0,R_{\mathrm{o}}) $  & $ \frac{1}{\sqrt{2}} $& $ \frac{1}{\sqrt{2}} $\\
		(1,3,3) & (2,3,3) & $ (0,0,R_{\mathrm{i}}) $  & $ (0,0,R_{\mathrm{o}}) $  & 1& 1 \\
		
		\bottomrule
	\end{tabular}
\end{table}
The covariance function in physical domain is defined as Eq.(\ref{Coveg2}) where $\sigma^2=1$, $b=1$ and $L=1$.

First, we decomposed the covariance function in the parametric space by using Algorithm \ref{ctt} with $tol=1.0000\times 10^{-5}$ and each $m_k=800$ in Algorithm \ref{ISE}, and computed the global relative errors $\varepsilon_{\mathrm{g}2}=3.3731\times 10^{-6}(N=1000)$ and ranks $\bm{r}_2=(4,21,60,21,4)$. Then, we computed the eigenpairs in each direction with $\varepsilon_{\mathrm{gf}2}$ = $1.3170\times 10^{-2}\ (N=1000)$.  Next, by setting $tol=10^{-4}$ and each $m_k=800$ in Algorithm \ref{ISE}, we computed the TT decomposition of the third-order cumulant function and obtained $\varepsilon_{\mathrm{g}3}=8.8170\times 10^{-4}(N=1000)$ and ranks $\bm{r}_3=(3,20,60,110,164,56,35,8)$. After HOSVD of $\bm{\breve{C}}_3$ ($tol_3=9.9990\times 10^{-1}$), the number of latent factors are reduced from 60 to 17 $ \varepsilon_{\mathrm{gf}3}=6.5547\times 10^{-3}, N=1000 $. The time cost (seconds) of each quantity is listed in Table \ref{eg3tb2}.
\begin{table}[htbp]
	\centering
	\caption{Time cost (seconds) of the quantities in Example 3}\label{eg3tb2}
	\begin{tabular}{*{6}{c}}
		\toprule
		&$ \tilde{C}_{k,\mathrm{TT}} $ & $\varepsilon_{\mathrm{g}k}$ & modes & $\varepsilon_{\mathrm{gf}k}$&total\\
		\midrule
		$ k =2$   & 202   & 1        & 2014   & 1 & 2218 \\ 
		$ k =3$   & 2821 & 11 & 31017 & 1 & 33850 \\
		\bottomrule
	\end{tabular}
\end{table}
 The results above indicate that the goal of matching the second and third-order cumulant functions within moderate time has been achieved.

\section{Conclusions}
In this work, we developed novel theoretical  and algorithm frameworks for extending the unidimensional  Karhunen-Lo{\`e}ve expansion to represent multidimensional random fields using higher-order cumulant functions. The algorithm framework can reveal the features of a random field in different directions in a \emph{highly automatic} way due to the marriage of a rank-revealing tensor train decomposition and the Chebfun paradigm of continuous computation. Differences from some existing methods were also discussed. Numerical experiments indicate that the proposed algorithm framework is able to overcome the current challenges of both computing the modes and representing the cumulants of the latent factors. The efficiency is moderately dependent on the parametric dimensionality, and is generally more than an order of magnitude higher than that of FEM. 

Finally, it is worth mentioning that the practical efficiency of an algorithm is also strongly affected by the computer languages. Computation of the third-order cumulant tensor of the latent factors contributes a major part of the time cost in each example. This phenomenon is mainly caused by the low efficiency of MATLAB on which the current Chebfun toolbox heavily depends. Hence, reconstructing the toolbox with a compiled language will also be beneficial for high-performance computing.

\section*{Acknowledgments}
This research was supported by the National Natural Science Foundation of China (Grant No.11572106) the National Key Research and Development Program of China (National Key Project No.2017YFC0703506), which are gratefully acknowledged by the authors. 

\bibliography{reference}
\end{document}